\documentclass[letterpaper]{amsart}
\pdfoutput=1
\newtheorem{theorem}{Theorem}[section]

\theoremstyle{definition}
\newtheorem{definition}[theorem]{Definition}

\theoremstyle{remark}

\usepackage{hyperref}
\usepackage[top=1.0in, left=1.0in, right=1.0in]{geometry}
\usepackage[all,cmtip]{xy}
\usepackage{xypic}
\usepackage{tikz,tikz-cd}
\usepackage{pgf}
\usetikzlibrary{positioning}
\usetikzlibrary{calc}
\usetikzlibrary{arrows}
\usetikzlibrary{shapes.multipart}
\usepackage{graphicx}
\usepackage{xcolor}
\definecolor{darkgreen}{rgb}{0.33, 0.55, 0.13}
\definecolor{darkpink}{rgb}{0.91, 0.33, 0.5}
\definecolor{electriccyan}{rgb}{0.0, 1.0, 1.0}
\definecolor{electricultramarine}{rgb}{0.25, 0.0, 1.0}
\definecolor{greenyellow}{rgb}{0.30, 0.82, 0.0}
\definecolor{lightbrown}{rgb}{0.71, 0.4, 0.11}
\definecolor{prune}{rgb}{0.45, 0.11, 0.11}
\usepackage{color}
\usepackage{amsmath}
\usepackage{amsfonts}
\usepackage{euler}
\usepackage{amssymb}
\usepackage{epsfig}
\usepackage{rotating}
\usepackage{graphics}
\newcommand{\be}{\begin{equation}}

\newcommand{\ee}{\end{equation}}

\newcommand{\dr}{{\mathrm{dir}}}

\newcommand{\dz}{\wedge}

\newcommand{\ba}{\begin{array}}

\newcommand{\ea}{\end{array}}

\newcommand{\beq}{\begin{eqnarray}}

\newcommand{\eeq}{\end{eqnarray}}

\newtheorem{lm}{lemma}

\newtheorem{thee}{theorem}

\newtheorem{proo}{proposition}

\newtheorem{co}{corollary}

\newtheorem{rem}{remark}

\newtheorem{deff}{definition}

\newcommand{\bd}{\begin{deff}}

\newcommand{\ed}{\end{deff}}

\newcommand{\bl}{\begin{lm}}

\newcommand{\el}{\end{lm}}

\newcommand{\bp}{\begin{proo}}

\newcommand{\ep}{\end{proo}}

\newcommand{\bt}{\begin{thee}}

\newcommand{\et}{\end{thee}}

\newcommand{\bc}{\begin{co}}

\newcommand{\ec}{\end{co}}

\newcommand{\brm}{\begin{rem}}

\newcommand{\erm}{\end{rem}}

\newcommand{\der}{{\rm d}}

\usepackage{t1enc}
\def\frak{\mathfrak}

\newcommand{\newc}{\newcommand}

\renewcommand{\exp}{\operatorname{exp}}

\let\ccdot\cdot
\def\cdot{\hbox to 2.5pt{\hss$\ccdot$\hss}}

\newc{\aR}{\mbox{\boldmath{$ R$}}}
\newc{\aS}{\mbox{\boldmath{$ S$}}}
\newc{\aT}{\mbox{\boldmath{$ T$}}}
\newc{\aW}{\mbox{\boldmath{$ W$}}}

\newc{\aK}{\mbox{\boldmath{$ K$}}}
\newc{\aL}{\mbox{\boldmath{$ L$}}}

\newcommand{\bbC}{\mathbb{C}}

\usepackage{amssymb}
\usepackage{amscd}

\newcommand{\hook}{\raisebox{-0.35ex}{\makebox[0.6em][r]
{\scriptsize $-$}}\hspace{-0.15em}\raisebox{0.25ex}{\makebox[0.4em][l]{\tiny
 $|$}}}

\newcommand{\bma}{\begin{pmatrix}}
\newcommand{\ema}{\end{pmatrix}}

\def\bbZ{{\mathbb{Z}}}

\newcommand{\X}{\mbox{\boldmath{$ X$}}}

\let\t=\tau
\let\m=\mu

\newc{\obstrn}[2]{B^{#1}_{#2}}

\newcommand{\rpl}                         
{\mbox{$
\begin{picture}(12.7,8)(-.5,-1)
\put(0,0.2){$+$}
\put(4.2,2.8){\oval(8,8)[r]}
\end{picture}$}}

\newcommand{\lpl}                         
{\mbox{$
\begin{picture}(12.7,8)(-.5,-1)
\put(2,0.2){$+$}
\put(6.2,2.8){\oval(8,8)[l]}
\end{picture}$}}

\usepackage{ifthen}

\newc{\tensor}[1]{#1}
\newc{\Mvariable}[1]{\mbox{#1}}
\newc{\down}[1]{{}_{#1}}
\newc{\up}[1]{{}^{#1}}

\newc{\JulyStrut}{\rule{0mm}{6mm}}
\newc{\midtenPan}{\mbox{\sf S}}
\newc{\midten}{\mbox{\sf T}}
\newc{\midtenEi}{\mbox{\sf U}}
\newc{\ATen}{\mbox{\sf E}}
\newc{\BTen}{\mbox{\sf F}}
\newc{\CTen}{\mbox{\sf G}}

\def\sideremark#1{\ifvmode\leavevmode\fi\vadjust{\vbox to0pt{\vss
 \hbox to 0pt{\hskip\hsize\hskip1em
 \vbox{\hsize3cm\tiny\raggedright\pretolerance10000
 \noindent #1\hfill}\hss}\vbox to8pt{\vfil}\vss}}}%

\newcommand{\Span}{\mathrm{Span}}

\numberwithin{equation}{section}

\newcounter{romenumi}
\newcommand{\labelromenumi}{(\roman{romenumi})}

\begin{document}
\newcommand{\bbS}{\mathbb{S}}
\newcommand{\bbR}{\mathbb{R}}
\newcommand{\bbK}{\mathbb{K}}
\newcommand{\sog}{\mathbf{SO}}
\newcommand{\spg}{\mathbf{Sp}}
\newcommand{\glg}{\mathbf{GL}}
\newcommand{\slg}{\mathbf{SL}}
\newcommand{\og}{\mathbf{O}}
\newcommand{\soa}{\frak{so}}
\newcommand{\spa}{\frak{sp}}
\newcommand{\gla}{\frak{gl}}
\newcommand{\sla}{\frak{sl}}
\newcommand{\sua}{\frak{su}}
\newcommand{\sug}{\mathbf{SU}}
\newcommand{\cspg}{\mathbf{CSp}}
\newcommand{\gat}{\tilde{\gamma}}
\newcommand{\Gat}{\tilde{\Gamma}}
\newcommand{\thet}{\tilde{\theta}}
\newcommand{\Thet}{\tilde{T}}
\newcommand{\rt}{\tilde{r}}
\newcommand{\st}{\sqrt{3}}
\newcommand{\kat}{\tilde{\kappa}}
\newcommand{\kz}{{K^{{~}^{\hskip-3.1mm\circ}}}}
\newcommand{\bv}{{\bf v}}
\newcommand{\di}{{\rm div}}
\newcommand{\curl}{{\rm curl}}
\newcommand{\cs}{(M,{\rm T}^{1,0})}
\newcommand{\tn}{{\mathcal N}}
\newcommand{\ten}{{\Upsilon}}
\title{A car as parabolic geometry}
\vskip 1.truecm
\author{C. Denson Hill} \address{Department of Mathematics, Stony Brook University, Stony Brook, NY 11794, USA}
\email{Dhill@math.stonybrook.edu}
\author{Pawe\l~ Nurowski} \address{Centrum Fizyki Teoretycznej,
Polska Akademia Nauk, Al. Lotnik\'ow 32/46, 02-668 Warszawa, Poland}
\email{nurowski@cft.edu.pl}
\thanks{Support: This work was supported by the Polish National Science Centre (NCN) via the grant number 2018/29/B/ST1/02583 and via the POLONEZ grant
2016/23/P/ST1/04148, which received funding from the European Union's Horizon
2020 research and innovation programme under the Marie Sk\l odowska-Curie grant
agreement No. 665778.}

\date{\today}
\begin{abstract}
We show that a car, viewed as a nonholonomic system, provides an example of a flat parabolic geometry of type $(\sog(2,3),P_{12})$, where $P_{12}$ is a Borel parabolic subgroup in $\sog(2,3)$. We discuss the relations of this geometry of a car with the geometry of circles in the plane (a low dimensional Lie sphere geometry), the geometry of 3-dimensional conformal Minkowski spacetime, the geometry of 3-rd order ODEs, projective contact geometry in three dimensions, and the corresponding twistor fibrations. We indicate how all these classical geometries can be interpreted in terms of the nonholonomic kinematics of a car.
\end{abstract}
\maketitle
\section{Car and Engel distribution}
\subsection{Configuration space and nonholonomic constraints}\label{intr}
In this note we look at a car from the point of view of an observer that is situated in space over the plane on which the car is moving. We idealize the car as an interval of length $\ell$ in the plane $\bbR^2$. The car has two pairs of wheels; we idealize them to be attached at both ends of the interval. The rear wheels are always parallel to the interval, whereas the front wheels can be rotated around the line vertical to the plane passing through the point of their attachment to the car. At every moment the direction of the front wheels can assume any angle with respect to the direction of the headlights of the car. To describe the position of the car we need \emph{four} numbers. One can define these four numbers in many ways; here we choose the setting depicted in the figure below:\\
\centerline{\includegraphics[height=6cm]{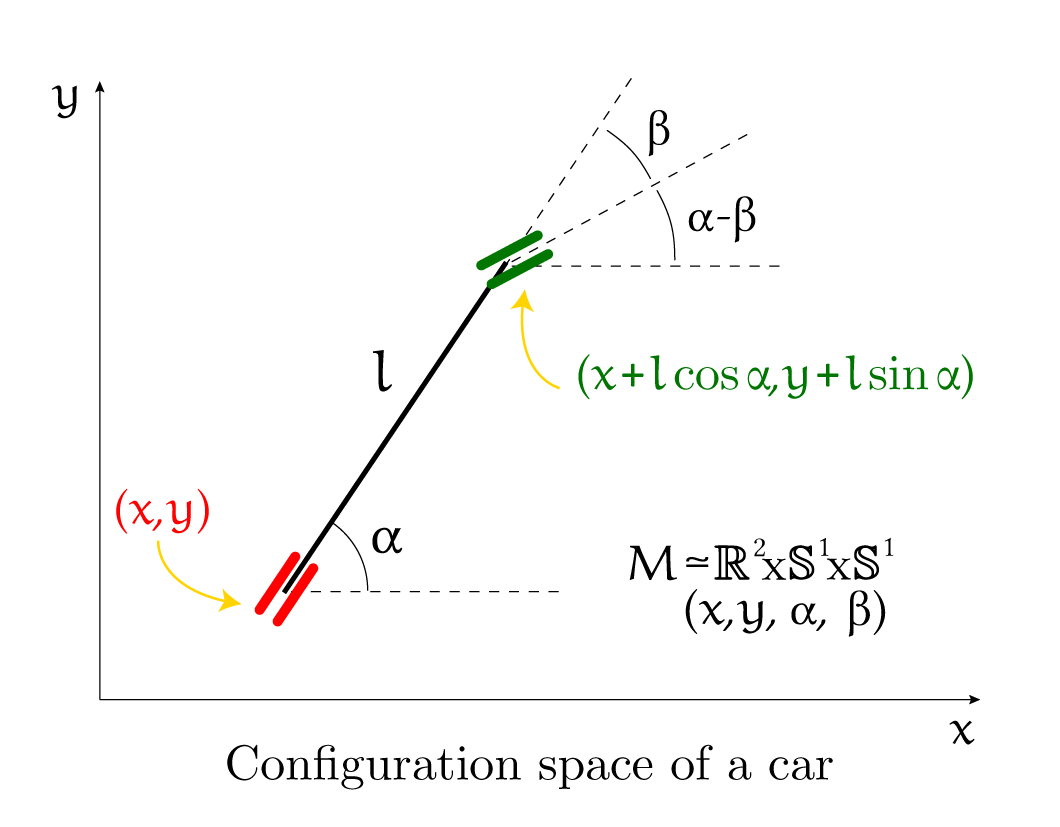}}
We introduce a Cartesian coordinate system in the plane so that the position of the rear wheels of the car has coordinates $(x,y)$. Then as a fixed line in the plane we choose the line $y=0$, and to keep track of the orientation of the chassis of the car we take the angle $\alpha$ that the interval representing the car forms with this line. The orientation of the front wheels is the angle $\beta$, between the direction defined by the front wheels and the direction of the interval representing the chassis of the car. As a result we have four numbers $(x,y,\alpha,\beta)$ describing uniquely the position of the car as it moves. Thus the \emph{configuration space} of the car is a 4-dimensional manifold $M$, locally diffeomorphic to $$\bbR^2\times\bbS^1\times\bbS^1=\{~(x,y,\alpha,\beta)~:~(x,y)\in\bbR^2;~\alpha,\beta\in\bbS^1~\}.$$

\subsection{Movement and the role of the tires}
When the car is moving it traverses a curve $q(t)=(x(t),y(t),$ $\alpha(t),\beta(t))$ in its configuration space $M$. The \emph{velocity} of the car at time $t$ is $\dot{q}(t)=(\dot{x}(t),\dot{y}(t),\dot{\alpha}(t),\dot{\beta}(t))$. It is a \emph{vector} from the tangent space $T_{q(t)}M$.

A safe car has \emph{tires}.  Their role is to prevent the car from \emph{skidding}. Our car will have \emph{perfect} tires. They impose \emph{nonholonomic constraints}. These are constraints on positions \emph{and velocities}, that can not be integrated to constraints on positions only. Indeed, what is expected from a properly behaving car is that its rear wheels, i.e. the point $(x,y)$ has its $(x,y)$-plane velocity parallel to the direction of the body of the car, and that the front wheels. i.e. the point $(x+\ell\cos\alpha,y+\ell\sin\alpha)$, has its $(x,y)$-plane velocity in the plane parallel to the orientation of the front wheels.
Thus, the movement of a car, represented by the curve $q(t)=(x(t),y(t),\alpha(t),\beta(t))\in M$, at every moment of time $t$, must satisfy
  $$\begin{aligned}
  &{ \color{red}\tfrac{\der}{\der t}(x,y)\quad||\quad(\cos\alpha,\sin\alpha)}\quad\quad\&\\& {\color{darkgreen}\tfrac{\der}{\der t}(x+\ell\cos\alpha,y+\ell\sin\alpha)\quad||\quad(\cos(\alpha-\beta),\sin(\alpha-\beta))},\end{aligned}$$
or, what is the same
$$\begin{aligned}
  &{ \color{red}\dot{x}\sin\alpha-\dot{y}\cos\alpha=0}\quad\quad\&\\& { \color{darkgreen}(\dot{x}-\ell\dot{\alpha}\sin\alpha)\,\sin(\alpha-\beta)-(\dot{y}+\ell\dot{\alpha}\cos\alpha)\,\cos(\alpha-\beta)=0}.\end{aligned}$$
 We emphasize that the above constraints are \emph{linear} in velocities.
Solving them we get the possible velocities as
$$\bma \dot{x}\\\dot{y}\\\dot{\alpha}\\\dot{\beta}\ema=A(t)\bma 0\\0\\0\\1\ema + B(t)\bma \ell \cos\alpha\cos\beta\\\ell\sin\alpha\cos\beta\\-\sin\beta\\0\ema.$$
where $\alpha=\alpha(t)$, $\beta=\beta(t)$, $A=A(t)$ and $B=B(t)$ are \emph{arbitrary} functions of time.
\subsection{Velocity distribution as an Engel distribution}
We can rephrase this by saying that at each point $q=(x,y,\alpha,\beta)^T$ in the tangent space $T_qM$, which is considered as the space of \emph{all} possible velocities, there is a \emph{distinguished} vector subspace\\
\centerline{\includegraphics[height=6cm]{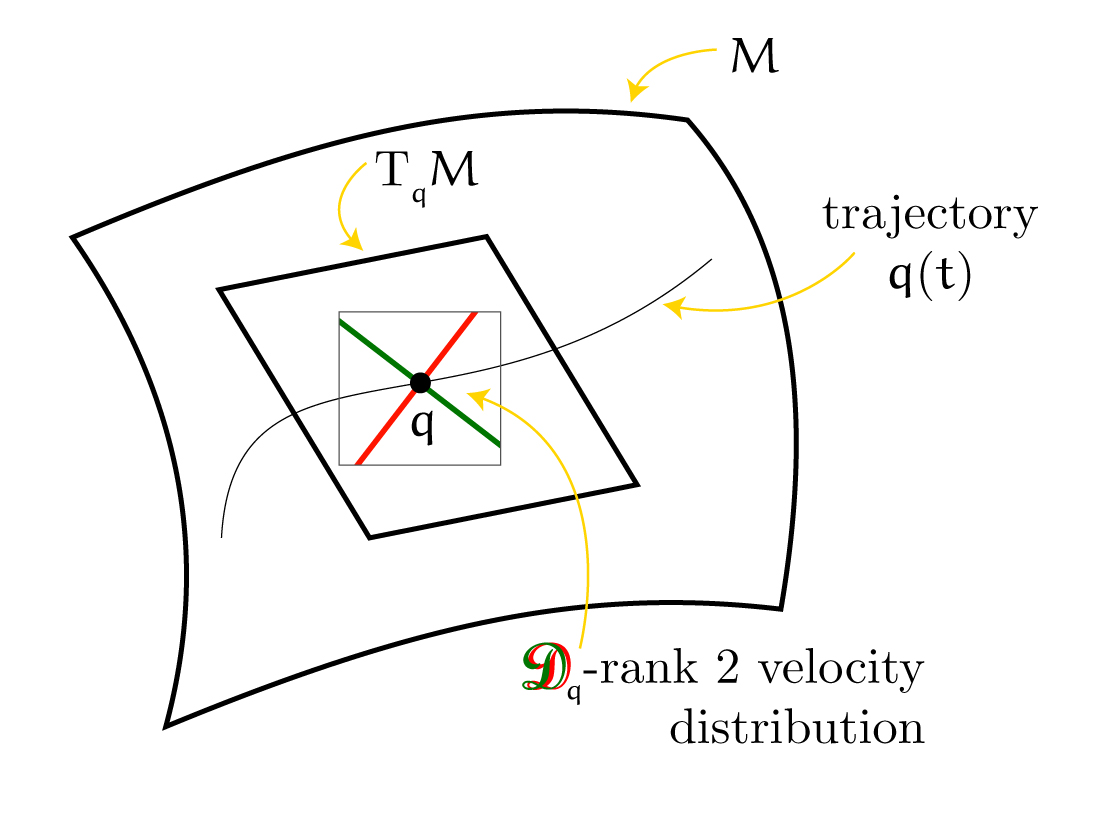}}
${\color{red}\mathcal D}\hspace{-0.31cm}{\color{darkgreen}\mathcal D}_q=\Span_\bbR({\color{darkgreen}X_3},{\color{red}X_4})$ spanned at each point $q\in M$ by the vectors tangent to the vector fields
 \be
   {\color{darkgreen}X_3=\partial_\beta}\quad\&\quad {\color{red}X_4=-\sin\beta\partial_\alpha+\ell\cos\beta(\cos\alpha\partial_x+\sin\alpha\partial_y)},\label{vd2}\ee
which is the space of \emph{admissible} velocities of the car at $q$. The car with perfect tires moves always along the curves $q(t)=(x(t),y(t),\alpha(t),\beta(t))^T$ such that its velocity $\dot{q}$ in the configuration space satisfies
$$\dot{q}=A {\color{darkgreen}X_3}+B {\color{red}X_4}.$$
The arbitrary functions $A=A(t)$ and $B=B(t)$ are called \emph{controls} of the car\footnote{Sometimes the vector fields $X_3$ and $X_4$ are also called controls.}.

Thus on $M$ there is a \emph{rank 2} distribution ${\color{red}\mathcal D}$\hspace{-0.31cm}${\color{darkgreen}\mathcal D}$ on $M$, describing the space of possible velocities, given by
\be{\color{red}\mathcal D}\hspace{-0.31cm}{\color{darkgreen}\mathcal D}=Span_{{\mathcal F}(M)}({ \color{darkgreen}X_3},{ \color{red}X_4}).\label{vd1}\ee
Therefore `the structure of a car with perfect tires' is \emph{up to now} $$(M,{\color{red}\mathcal D}\hspace{-0.3cm}{\color{darkgreen}\mathcal D}),$$
i.e. a 4-manifold $M$ with a rank 2 distribution $(M,{\color{red}\mathcal D}\hspace{-0.31cm}{\color{darkgreen}\mathcal D})$. 

Now the fundamental question is: \emph{Is} ${\color{red}\mathcal D}\hspace{-0.3cm}{\color{darkgreen}\mathcal D}$ \emph{integrable?}

The answer is: Obviously \emph{not}, since everybody knows that a car can be driven from any position in its configuration space to any other position (Chow-Raszewski theorem). One can also convince oneself about that by calculating the commutators of ${\color{darkgreen}X_3}$ and ${\color{red}X_4}$. We have:
   \be\begin{aligned}
  {}[{\color{darkgreen}X_3},{\color{red}X_4}]&=-\cos\beta\partial_\alpha-\ell\sin\beta(\sin\alpha\partial_y+\cos\alpha\partial_x):=X_2\\
    [{\color{red}X_4},X_2]&=\ell(\cos\alpha\partial_y-\sin\alpha\partial_x):=X_1,\end{aligned}\label{vd4}
\ee
and it is easy to check that
$$X_1\wedge X_2\wedge {\color{darkgreen}X_3}\wedge {\color{red}X_4}=\ell^2\partial_x\wedge\partial_y\wedge\partial_\alpha\wedge\partial_\beta\neq 0.$$
This shows that taking successive commutators of the vectors from the car distribution ${\color{red}\mathcal D}\hspace{-0.3cm}{\color{darkgreen}\mathcal D}$ we quickly (in two steps!) produce the entire tangent bundle to $M$. This, by the Chow-Raszewski theorem, is a well know condition for curves tangent to the distribution to be capable reaching any point of the configuration space from any other point.

We summarize this by defining three distributions ${\mathcal D}_{-1}$, ${\mathcal D}_{-2}$ and ${\mathcal D}_{-3}$ on $M$ as in the table below:
$$
  \begin{matrix}
    &&\mathrm{rank}\\
    {\mathcal D}_{-1}:={\color{red} \mathcal D}\hspace{-0.31cm}{\color{darkgreen} \mathcal D} & \mathrm{Span}({\color{red}X_4},{\color{darkgreen}X_3}) & {\color{blue}2}\\
    {\mathcal D}_{-2}:=[{\mathcal D}_{-1},{\mathcal D}_{-1}]& \mathrm{Span}({\color{red}X_4},{\color{darkgreen}X_3},X_2)& {\color{blue}3}\\
    {\mathcal D}_{-3}:=[{\mathcal D}_{-1},{\mathcal D}_{-2}]& \mathrm{Span}({\color{red}X_4},{\color{darkgreen}X_3},X_2,X_1)=\mathrm{T}M& {\color{blue}4}
\end{matrix}
$$
Thus given the so far defined structure of the car $(M,{\color{red}\mathcal D}\hspace{-0.3cm}{\color{darkgreen}\mathcal D})$, we have a filtration ${\mathcal D}_{-1}\subset {\mathcal D}_{-2}\subset {\mathcal D}_{-3}=\mathrm{T}M$ of distributions with the \emph{constant growth vector} ${\color{blue}(2,3,4)}$. These collective properties of  the car distribution ${\color{red}\mathcal D}\hspace{-0.3cm}{\color{darkgreen}\mathcal D}$ make it an \emph{Engel distribution}. Here we recall that an abstract \emph{Engel distribution} is a rank 2 distribution on a 4-manifold such that its derived flag of distributions ${\mathcal D}_{-1}={\mathcal D}$, ${\mathcal D}_{-2}:=[{\mathcal D}_{-1},{\mathcal D}_{-1}]$ and ${\mathcal D}_{-3}:=[{\mathcal D}_{-1},{\mathcal D}_{-2}]$ has respective \emph{constant} ranks $2,3$ and $4$. 
\subsection{Equivalence of Engel distributions}
Our discussion so far shows that the geometric structure associated with a car is $(M,{\color{red}\mathcal D}\hspace{-0.31cm}{\color{darkgreen}\mathcal D})$ with ${\color{red}\mathcal D}\hspace{-0.31cm}{\color{darkgreen}\mathcal D}$ being an Engel distribution on a manifold $M$.

A newcomer to this subject has an immediate question: are there nonequivalent Engel distributions? To answer this we need the notion of \emph{equivalence} of distributions.

We say that two distributions ${\mathcal D}$ and $\bar{\mathcal D}$ of the same rank on manifolds $M$ and $\bar{M}$ of the same dimension are \emph{(locally) equivalent} iff there exists a (local) diffeomorphism $\phi:M\to \bar{M}$ such that $\phi_*{\mathcal D}=\bar{\mathcal D}$. (Local) self-equivalence maps $\phi:M\to M$, i.e. maps such that $\phi_*{\mathcal D}={\mathcal D}$  are called (local) \emph{symmetries} of $\mathcal D$. They form a \emph{group of (local) symmetries of } $\mathcal D$. This notion has its infinitesimal version: we say that a vector field $X$ on $M$ is an \emph{infinitesimal symmetry} of $\mathcal D$ if and only if ${\mathcal L}_X{\mathcal D}\subset{\mathcal D}$. Since the commutator $[X,Y]$ of two infinitesimal symmetries $X$ and $Y$ is also an infinitesimal symmetry, this leads to the notion of the \emph{Lie algebra} $\mathfrak{g}_{\mathcal D}$ \emph{of infinitesimal symmetries} of $\mathcal D$.

Now, one convinces herself that the distribution  
$${\mathcal D}_E=({\color{darkgreen}\partial_q},{\color{red}\partial_x+p\partial_y+q\partial_p})$$ defined on an open set of $\bbR^4$ parametrized by $(x,y,p,q)$ is an Engel distribution. We have the following classical theorem due to Friedrich Engel.

{\bf Theorem} Every Engel distribution is locally equivalent to the distribution ${\mathcal D}_E$.

  One may say that we are in trouble: Since the car structure $(M,{\color{red}\mathcal D}\hspace{-0.31cm}{\color{darkgreen}\mathcal D})$ is a structure of an Engel distribution, there is no geometry associated to the car. The wrong argument in this kind of criticisim is that an Engel distribution is not the only structure that a car with perfect tires has. It turns out that the geometry associated with a car is more subtle than just the geometry of an Engel distribution. The car features equip its Engel distribution with an additional structure.

  \section{Car and Engel distribution with a split}
   \subsection{Two distinguished directions}
To see this consider the vector field:
    ${ \color{red}X_4=-\sin\beta\partial_\alpha+\ell\cos\beta(\cos\alpha\partial_x+}$ ${ \color{red}\sin\alpha\partial_y)}$. When $\beta=0$ it becomes
${\color{red}X_4=\ell(\cos\alpha\partial_x+\sin\alpha\partial_y)}$ and if the car chooses this direction of its velocity it makes a simple movement by going along a straight line in the direction $(\cos\alpha,\sin\alpha)$ in the $(x,y)$ plane. On the other hand,  if the car chooses its velocity in the direction of the vector fild ${\color{darkgreen}X_3=\partial_\beta}$, then although it does move in the configuration space, it does not perform any movement in the physical $(x,y)$ plane, merely rotating the steering wheel/front wheels with the engine at idle.

Cars owners/producers perfectly know and \emph{make use} of these two particular vector fields $({\color{darkgreen}X_3},{\color{red}X_4})$ in the distribution ${\color{red}\mathcal D}\hspace{-0.31cm}{\color{darkgreen}\mathcal D}$. In particular, car owners alternate using these two vector fields, each separately at proper instants/intervals of time, in parallel parking.

Indeed, if one wants to park a car one first approaches the parking spot by having its velocity aligned with ${\color{red}X_4}$ vector field with $\beta=0$. Then the car stops and rotates its front wheels towards the sidewalk passing from $\beta=0$ to $\beta=\beta_0$=const. This is done by aligning its velocity with the vector field ${\color{darkgreen}X_3}$. After this, the car velocity again becomes aligned with ${\color{red}X_4}$, which now has $\beta=\beta_0$=const, so that the car goes backwards towards the sidewalk.\\
\centerline{\includegraphics[height=6cm]{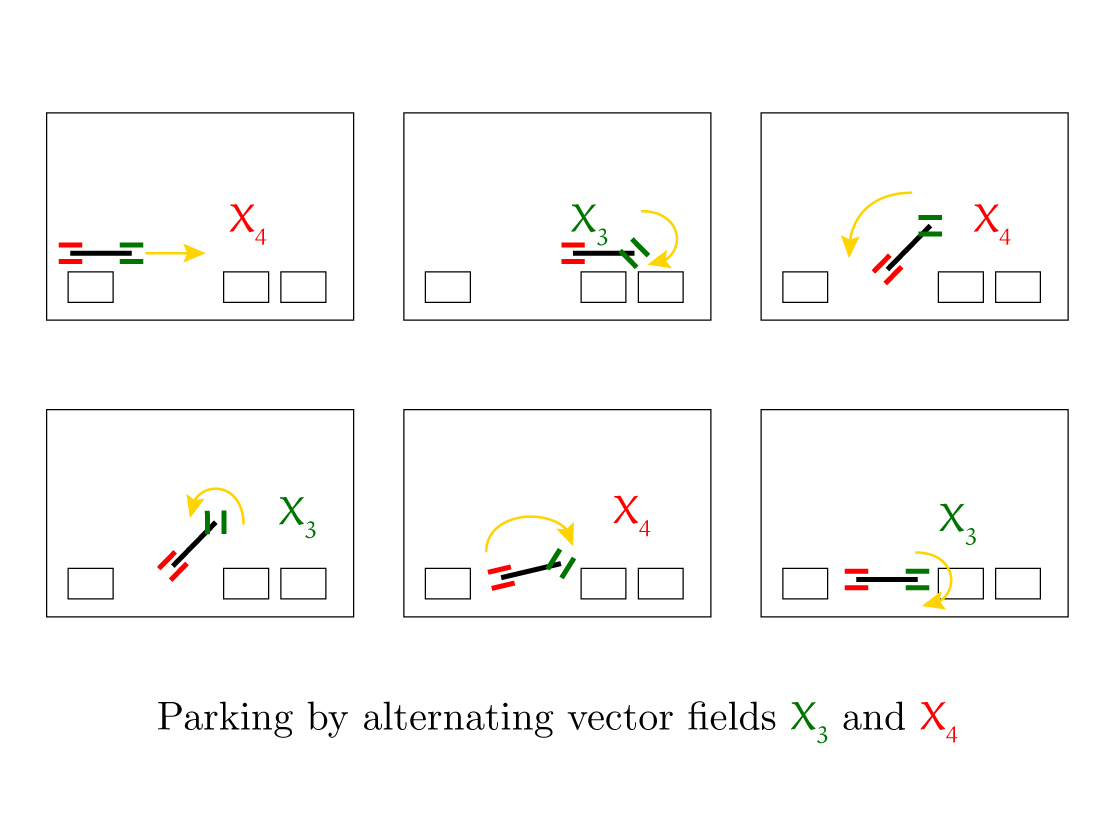}}
When the rear wheels are close to the sidewalk the car stops again, and aligns its velocity with ${\color{darkgreen}X_3}$, going back from $\beta=\beta_0$ to $\beta=-\beta_0$. Again applying backwards ${\color{red} X_4}$ with this constant $\beta=-\beta_0$ enables the driver to orient the rear wheels parallely to the sidewalk. If this happens, the car stops and applies ${\color{darkgreen}X_3}$ to make $\beta=0$ again. Finally the car aligns its velocity with ${\color{red}X_4}$ having $\beta=0$ to move parallely to the sidewalk and to take the midlle position between the two cars before and after it. 

Thus the car's distribution ${\color{red}\mathcal D}\hspace{-0.31cm}{\color{darkgreen}\mathcal D}$ has an additional structure, which is its split
$${\color{red}\mathcal D}\hspace{-0.31cm}{\color{darkgreen}\mathcal D}={\color{darkgreen}{\mathcal D}_w}\oplus {\color{red}{\mathcal D}_g},$$
onto rank one subdistributions
$${\color{darkgreen}{\mathcal D}_w}=\Span_{{\mathcal F}(M)}({\color{darkgreen}X_3})\quad\mathrm{and}\quad {\color{red}{\mathcal D}_g}=\Span_{{\mathcal F}(M)}({\color{red}X_4}).$$
    These subdistributions have a clear physical meaning:

The distribution  ${\color{darkgreen}{\mathcal D}_w}$ as spanned by ${\color{darkgreen}X_3=\partial_\beta}$, is responsible for the steering wheel control, and will be called the \emph{steering wheel space}; on the other hand  the distribution ${\color{red}{\mathcal D_g}}$, as spanned by the generator of the forward-backward movement ${ \color{red}X_4=-\sin\beta\partial_\alpha+\ell\cos\beta(\cos\alpha\partial_x+\sin\alpha\partial_y)}$ will be called the \emph{gas space}.
      
      This results in the statement that the car structure is actually $(M,{\color{red}\mathcal D}\hspace{-0.31cm}{\color{darkgreen}\mathcal D}={\color{darkgreen}{\mathcal D}_w}\oplus {\color{red}{\mathcal D}_g})$, with ${\color{red}\mathcal D}\hspace{-0.31cm}{\color{darkgreen}\mathcal D}$ being an Engel distribution with a \emph{split} ${\color{red}\mathcal D}\hspace{-0.31cm}{\color{darkgreen}\mathcal D}={\color{darkgreen}{\mathcal D}_w}\oplus {\color{red}{\mathcal D}_g}$ onto rank one, steering wheel and gas, subdistributions. So considering a car's geometry more thoroughly we land in a realm of the subtle geometry of Engel distributions \emph{with a split}!
  
\subsection{New geometry: Engel distributions with a split}
Thus we ultimately established that the geometry of a car with perfect tires, is given by a structure $(M,{\color{red}\mathcal D}\hspace{-0.31cm}{\color{darkgreen}\mathcal D}={\color{darkgreen}{\mathcal D}_w}\oplus {\color{red}{\mathcal D}_g})$, where ${\color{red}\mathcal D}\hspace{-0.31cm}{\color{darkgreen}\mathcal D}$ is an Engel distribution \emph{with a (car's) split} ${\color{red}\mathcal D}\hspace{-0.31cm}{\color{darkgreen}\mathcal D}={\color{darkgreen}{\mathcal D}_w}\oplus {\color{red}{\mathcal D}_g}.$

Abstractly, irrespectively of car's considerations, let us consider a geometry in the form $(M, {\mathcal D}={\mathcal D}_1\oplus {\mathcal D}_2)$, where dim$M$=4, $\mathcal D$ is an Engel distribition on $M$, and both subdistributions ${\mathcal D}_1$ and ${\mathcal D}_2$ in $\mathcal D$ have rank one. Let us call this an \emph{Engel structure with a split}.

Such structures have their own equivalence problem, related to the following definitions:

Two Engel structures with a split $(M, {\mathcal D}={\mathcal D}_1\oplus {\mathcal D}_2)$ and $(\bar{M}, \bar{{\mathcal D}}=\bar{{\mathcal D}_1}\oplus \bar{{\mathcal D}_2})$ are \emph{(locally) equivalent} if and only if there exists a (local) diffeomorphism $\phi:M\to \bar{M}$ such that $\phi_*{\mathcal D}_1=\bar{{\mathcal D}_1}$ and $\phi_*{\mathcal D}_2=\bar{{\mathcal D}_2}$.
Infinitesimally, we consider vector fields $S$ on $M$ such that ${\mathcal L}_S{\mathcal D}_1\subset{\mathcal D}_1$ and ${\mathcal L}_S{\mathcal D}_2\subset{\mathcal D}_2$, and we call such vector fields \emph{infinitesimal symmetries} of $(M, {\mathcal D}={\mathcal D}_1\oplus {\mathcal D}_2)$. This, as usual, leads to a notion of the \emph{Lie algebra} $\mathfrak{g}_{\mathcal D}$ \emph{of infinitesimal symmetries of} an Engel structure $(M, {\mathcal D}={\mathcal D}_1\oplus {\mathcal D}_2)$ with a split, as the Lie algebra of the vectors fields $S$ as above.

We can now ask about the Lie algebra of infinitesimal symmetries of the Engel structure with a split  $(M,{\color{red}{\mathcal D}}\hspace{-0.31cm}{\color{darkgreen}{\mathcal D}}={\color{darkgreen}{\mathcal D}_w}\oplus {\color{red}{\mathcal D}_g})$ of a car. In this case we have ${\mathcal D}_1={\color{darkgreen}{\mathcal D}_w}$ and ${\mathcal D}_2={\color{red}{\mathcal D}_g}$. As an answer we get a bit surprising result as below: 
\begin{theorem}
  Consider the car structure $(M,{\color{red}{\mathcal D}}\hspace{-0.31cm}{\color{darkgreen}{\mathcal D}})$ consisting of its velocity distribution ${\color{red}{\mathcal D}}\hspace{-0.31cm}{\color{darkgreen}{\mathcal D}}$ and the split of ${\color{red}{\mathcal D}}\hspace{-0.31cm}{\color{darkgreen}{\mathcal D}}$ onto rank 1 distributions ${\color{red}{\mathcal D}}\hspace{-0.31cm}{\color{darkgreen}{\mathcal D}}={\color{darkgreen}{\mathcal D_w}}\oplus {\color{red}{\mathcal D_g}}$ with 
$
{\color{darkgreen}{\mathcal D}_w=\mathrm{Span}(\partial_\beta)}, \,\,  {\color{red}{\mathcal D}_g=\mathrm{Span}(-\sin\beta\partial_\alpha+\ell\cos\beta(\cos\alpha\partial_x+\sin\alpha\partial_y)}.$\\
The Lie algebra of infinitesimal symmetries of this Engel structure with a split is 10-dimensional, with the following generators
 $$\begin{aligned}
    S_1&=\partial_x\\
    S_2&=\partial_y\\
    S_3&=x\partial_y-y\partial_x+\partial_\alpha\\
    S_4&=\ell(\sin\alpha\partial_x-\cos\alpha\partial_y)+\sin^2\beta\partial_\beta\\
    S_5&=x\partial_x+y\partial_y-\sin\beta\cos\beta\partial_\beta\\
    S_6&=(x^2-y^2)\partial_x+2xy\partial_y+2y\partial_\alpha-2\cos\beta\Big(\ell\cos\beta\sin\alpha+x\sin\beta\Big)\partial_\beta\\
    S_7&=\ell\Big(x(\sin\alpha\partial_x-\cos\alpha\partial_y)-\cos\alpha\partial_\alpha\Big)+\sin\beta\Big(\ell\cos\beta\sin\alpha+x\sin\beta\Big)\partial_\beta\\
    S_8&=\ell\Big(y(\sin\alpha\partial_x-\cos\alpha\partial_y)-\sin\alpha\partial_\alpha\Big)-\sin\beta\Big(\ell\cos\beta\cos\alpha-y\sin\beta\Big)\partial_\beta\\
    S_9&=2xy\partial_x+(y^2-x^2)\partial_y-2x\partial_\alpha+2\cos\beta\Big(\ell\cos\beta\cos\alpha-y\sin\beta\Big)\partial_\beta\\
    S_{10}&=\ell(x^2+y^2)\Big(\sin\alpha\partial_x-\cos\alpha\partial_y\Big)-2\ell\Big(x\cos\alpha+y\sin\alpha\Big)\partial_\alpha+\\&\Big(2\ell\sin\beta\cos\beta\big(x\sin\alpha-y\cos\alpha\big)+\sin^2\beta(x^2+y^2)+2\ell^2\cos^2\beta\Big)\partial_\beta
\end{aligned}$$
It is isomorphic to the simple real Lie algebra $\mathfrak{so}(2,3)=\mathfrak{sp}(2,\bbR)$.  Moreover, there are plenty of locally nonequivalent Engel distributions with a split, but the split ${\color{red}{\mathcal D}}\hspace{-0.31cm}{\color{darkgreen}{\mathcal D}}={\color{darkgreen}{\mathcal D}_w}\oplus {\color{red}{\mathcal D}_g}$ on the (Engel) car distribution used  by car owners and provided by cars' producers is \emph{the most symmetric}.
\end{theorem}
The fact that there are many locally nonequivalent Engel structures with a split is not surprising at all. What is surprising here, is that the split on the Engel distribution provided by the `steering-wheel--gas' control of a car is the most symmetric. Moreover, the appearence of a \emph{simple} Lie algebra $\mathfrak{so}(2,3)=\mathfrak{sp}(2,\bbR)$ as the full algebra of symmetries of car's ${\color{red}{\mathcal D}}\hspace{-0.31cm}{\color{darkgreen}{\mathcal D}}={\color{darkgreen}{\mathcal D_w}}\oplus {\color{red}{\mathcal D_g}}$ is also striking. Especially that $\mathfrak{so}(2,3)$ is the Lie algebra of the group of conformal symmetries of 3-dimensional \emph{Minkowski space}. How on earth Minkowski space can be related to a car?
  \section{Explaining the $\mathfrak{so}(2,3)=\mathfrak{sp}(2,\bbR)$ symmetry}
 \subsection{A  double fibration}\label{dfib}
  Consider integral curves of the two distinguished directions ${\color{darkgreen}X_3}$ and ${\color{red}X_4}$ defined by the split in the car's distribution ${\color{red}{\mathcal D}}\hspace{-0.31cm}{\color{darkgreen}{\mathcal D}}$. Let us call the integral curves of ${\color{darkgreen}X_3}$ by ${\color{darkgreen}q_3}$ and the integral curves of ${\color{red}X_4}$ by ${\color{red}q_4}$ respectively. They define two \emph{foliations} of $M$, the first having ${\color{darkgreen}q_3}$ as the leaves, and the second consisting of leaves given by ${\color{red}q_4}$. Passing to the space of leaves of these two foliations, which we denote by ${\color{darkgreen}P}$ and by ${\color{red}Q}$, respectively, we get a double fibration\\
  \centerline{\includegraphics[height=5cm]{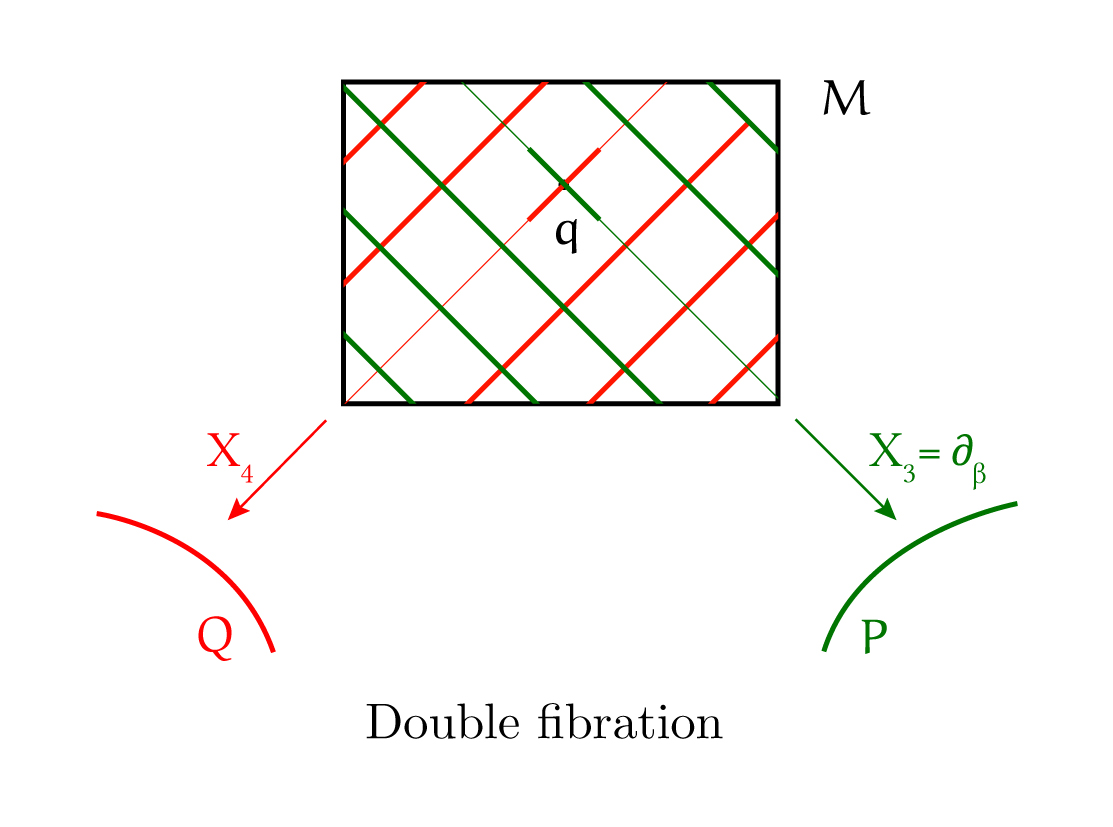}}
  with the 4-dimensional configuration space $M$ of a car on top, and the two 3-dimensional spaces  ${\color{darkgreen}P}$ and ${\color{red}Q}$ at the bottom.

  We will now analyze the geometry of each of the base spaces of this fibration, devoting a subsection to each of them. 
  \subsection{Conformal structure on ${\color{red}Q}$}\label{sec41}
  Points of ${\color{red}Q}$ are just the integral curves of ${\color{red}X_4}$.  What are these curves in $M$?  In an appropriate parametrization they are:
\be {\color{red}q_4(t)=\bma 2\ell\cot\beta_0\cos(\alpha_0-\tfrac12 t\sin\beta_0)\sin(\tfrac12t\sin\beta_0)+x_0\\
2\ell\cot\beta_0\sin(\alpha_0-\tfrac12 t\sin\beta_0)\sin(\tfrac12 t\sin\beta_0)+y_0\\
-t\sin\beta_0+\alpha_0\\\beta_0\ema}\quad\mathrm{when}\quad\beta_0\neq 0,\label{q4a}\ee
or
\be {\color{red}q_4(t)=\bma t\ell\cos\alpha_0+x_0\\
t\ell\sin\alpha_0+y_0\\
\alpha_0\\0\ema}\quad\mathrm{when}\quad\beta_0=0.\label{q4b}\ee
Here $(x_0,y_0,\alpha_0,\beta_0)$ are constants, corresponding to the position of the car at $t=0$.

These curves ${\color{red}q_4(t)}$  correspond to the movement of the car, when the $\beta$ angle is fixed. Thus in the configuration space $M$, they are \emph{helices} $(x(t),y(t),\alpha(t))$ in the 3-dimensional space $\beta=\beta_0$=const, parametrized by $(x,y,\alpha)$. The axi of these helices are given by $(x_0+\ell\cot\beta_0\sin\alpha_0,y_0-\ell\cot(\beta_0)\cos\alpha_0,t)$, their radii are $R=\ell\cot\beta_0$ and their pitch is $2\pi$, for each choice of initial conditions $(x_0,y_0,\alpha_0)$.

 \centerline{\includegraphics[height=5cm]{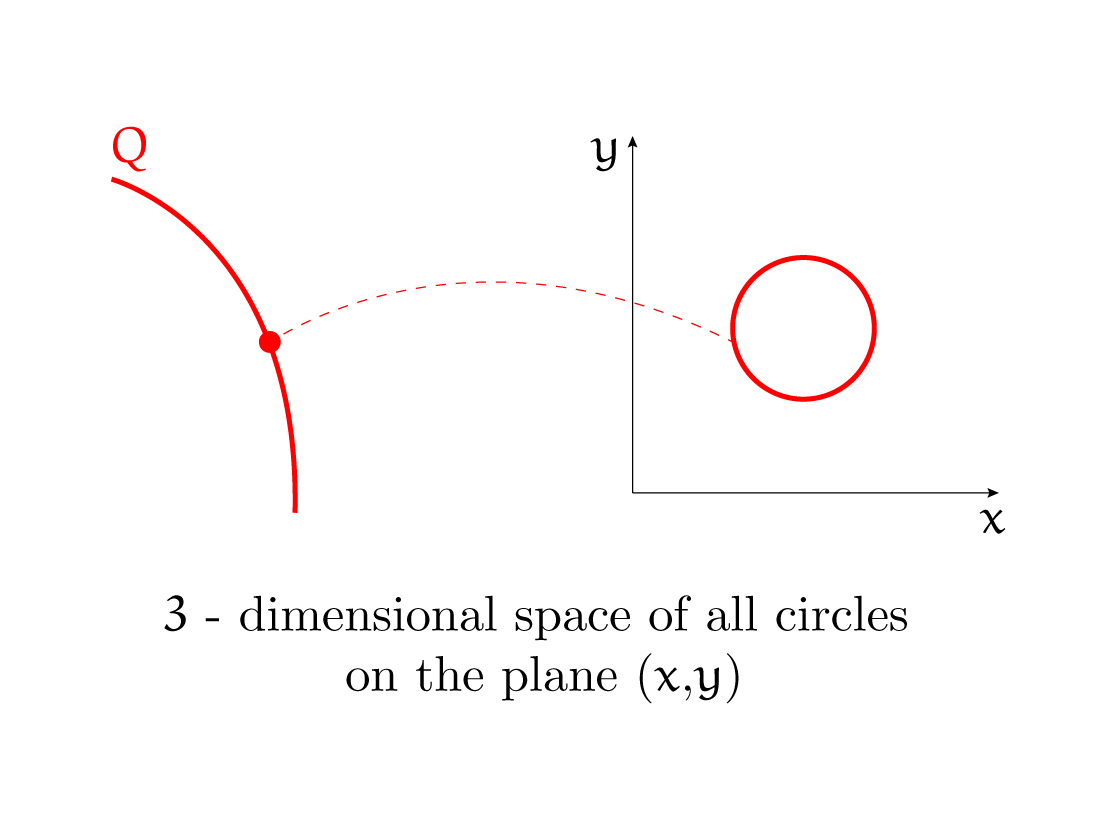}}

In the physical 2-dimensional space $(x,y)$, where the car is physically moving, these curves are either points (when $\beta_0=\pm \pi/2$), or circles (when $0<|\beta_0|<\pi/2$), or straight lines (when $\beta_0=0$). This corresponds to the simple fact that if one sets the steering wheel in a given position, or what is the same keeps the constant angle $\beta=\beta_0$ between the front wheels and the axis of the chasis of the car, the \emph{rear} wheels of the car will go on a straight line if $\beta=0$, will go on circles if $0<|\beta|<\pi/2$, or will stay at a given point $(x_0,y_0)$ if the front wheels are perpendicular to the axis of the car. It is important to note that, by setting the initial conditions $(x_0,y_0,\alpha_0,\beta_0)$ properly, one can obtain any point, line or a circle in the plane $(x,y)$, as a trajectory of a physical movement of the car in the plane $(x,y)$.

Thus there is a one-to one correspondence between the points ${\color{red}q}$ of the 3-dimensional space ${\color{red}Q}$ of the integral curves of the vector field ${\color{red}X_4}$ (the helices at each plane $\beta=\beta_0$ in $M$) and the 3-dimensional \emph{space} ${\color{red}{\bf Q}}$ \emph{of all points, circles and lines in} $\bbR^2$ coordinatized by $(x,y)$.   

\subsubsection{Geometry of \emph{oriented} circles on the plane}

Since two circles on the plane can be disjoint, or can intersect, or be tangent, and since these relations between any two circles are invariant with respect to diffeomorphisms of the plane, they should be used to further determine the geometry of the space ${\color{red}{\bf Q}}$ and in turn the geometry of the leaf space ${\color{red}Q}$.


The geometry of circles on the plane is a classical subject first considered by S. Lie (see e.g. \cite{helgason}). Consider a set ${\color{red}{\bf Q}}$ of all objects in the plane whose coordinates $(x,y)$ satisfy $$x^2+y^2-2ax-2by+c=0,$$
with some real constants $a,b,c$. Introducing
$$R^2=a^2+b^2-c,$$
and projective coordinates $[\xi:\eta:\zeta:\mu:\nu]$ in $\bbR P^4$ via
\be a=\frac{\xi}{\nu},\quad b=\frac{\eta}{\nu},\quad c=\frac{\mu}{\nu},\quad R=\frac{\zeta}{\nu},\label{abc}\ee
we see that ${\color{red}{\bf Q}}$ is a \emph{projective quadric}
\be{\color{red}{\bf Q}}~=~\{~\bbR P^4\ni[\xi:\eta:\zeta:\mu:\nu]~:~~\xi^2+\eta^2-\zeta^2-\mu\nu=0~\}\label{abcd}\ee
in $\bbR P^4$. 
The objects (the points) of this set are stratified as follows. Generically they form the set ${\color{red}{\bf Q}_c}$ of (all) circles in the plane; this occurs when $\xi^2+\eta^2-\mu\nu>0$. When the radius $R$ is infinite, i.e. when $\nu=0$, the objects belong to ${\color{red}{\bf Q}_\ell}$, the set of (all) lines in the plane; finally, when $\zeta=0$, the objects belong to ${\color{red}{\bf Q}_p}$, the set of (all) points on the plane. Thus we have
$${\color{red}{\bf Q}={\bf Q}_c\sqcup{\bf Q}_\ell\sqcup{\bf Q}_p},$$
i.e. ${\color{red}{\bf Q}}$ is the set of \emph{all circles}, \emph{lines} and \emph{points} on the plane. In addition we easily see that the \emph{three dimensional} set ${\color{red}{\bf Q}}$, as a \emph{null} projective quadric in $\bbR P^4$, acquires a natural \emph{conformal Lorentzian structure} $[g]$, coming from the quadratic form
\be
Q(\xi,\eta,\zeta,\mu,\nu)=\xi^2+\eta^2-\zeta^2-\mu\nu\label{qform}\ee
in $\bbR^5$. 

It is important to notice that by considering $R$ as in formula \eqref{abc} we \emph{doubled} the number of circles in the plane. This is because, depending on the sign of $\zeta\nu$, the radius $R$ of the circle may be positive or negative. This has an obvious interpretation: the space $\color{red}{\bf Q}$ consists of all \emph{oriented} circles/lines. We adapt the convention that a circle/line $(x-a)^2+(y-b)^2=R^2$ is oriented \emph{counterclockwise} iff $R>0$, and it is oriented \emph{clockwise} iff $R<0$.    

\centerline{\includegraphics[height=5cm]{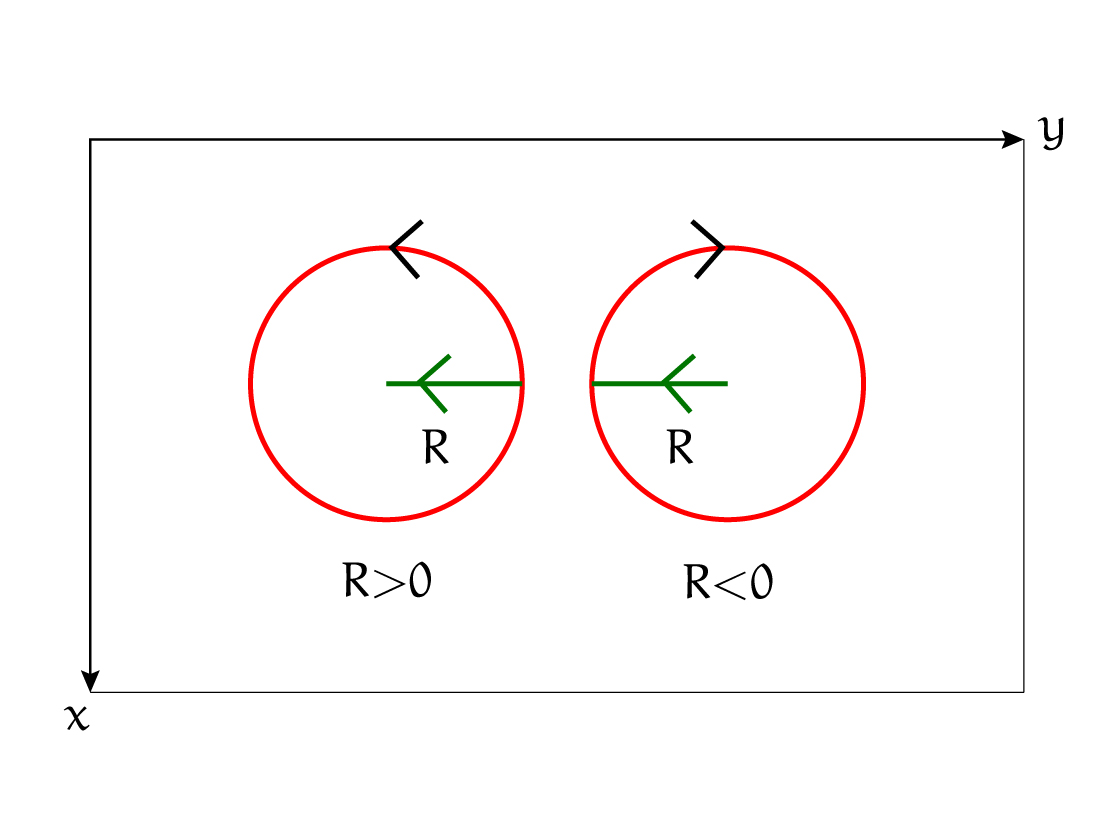}  \includegraphics[height=5cm]{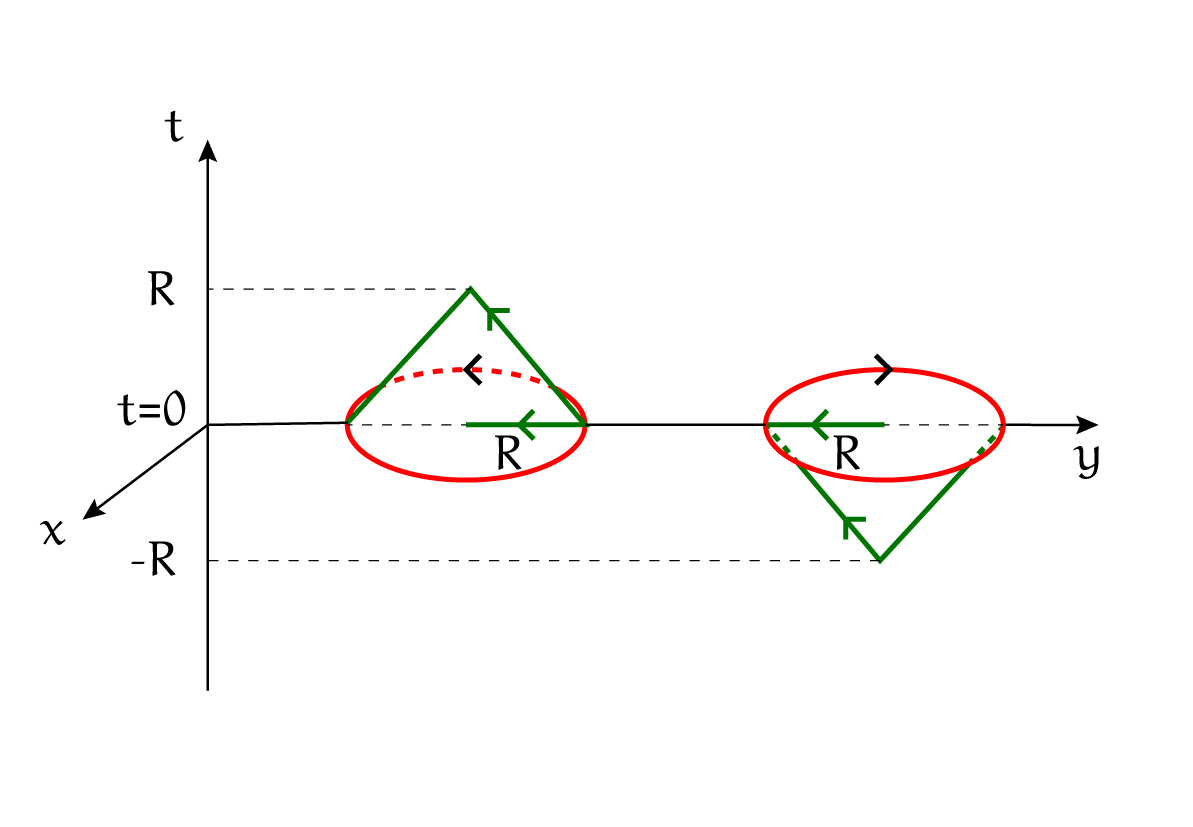}}

Lie has shown that the conformal structure in ${\color{red}{\bf Q}}$, whose points are generically oriented circles in the plane, is identical with the structure defined by the \emph{incidence} relation between the circles: \emph{two circles} from ${\color{red}{\bf Q}}$ \emph{are incident if and only if they are tangent to each other in such a way that their orientations coincide when one of the circles is inside the other and are opposite when they are external to each other}.\\
\centerline{\includegraphics[height=7cm]{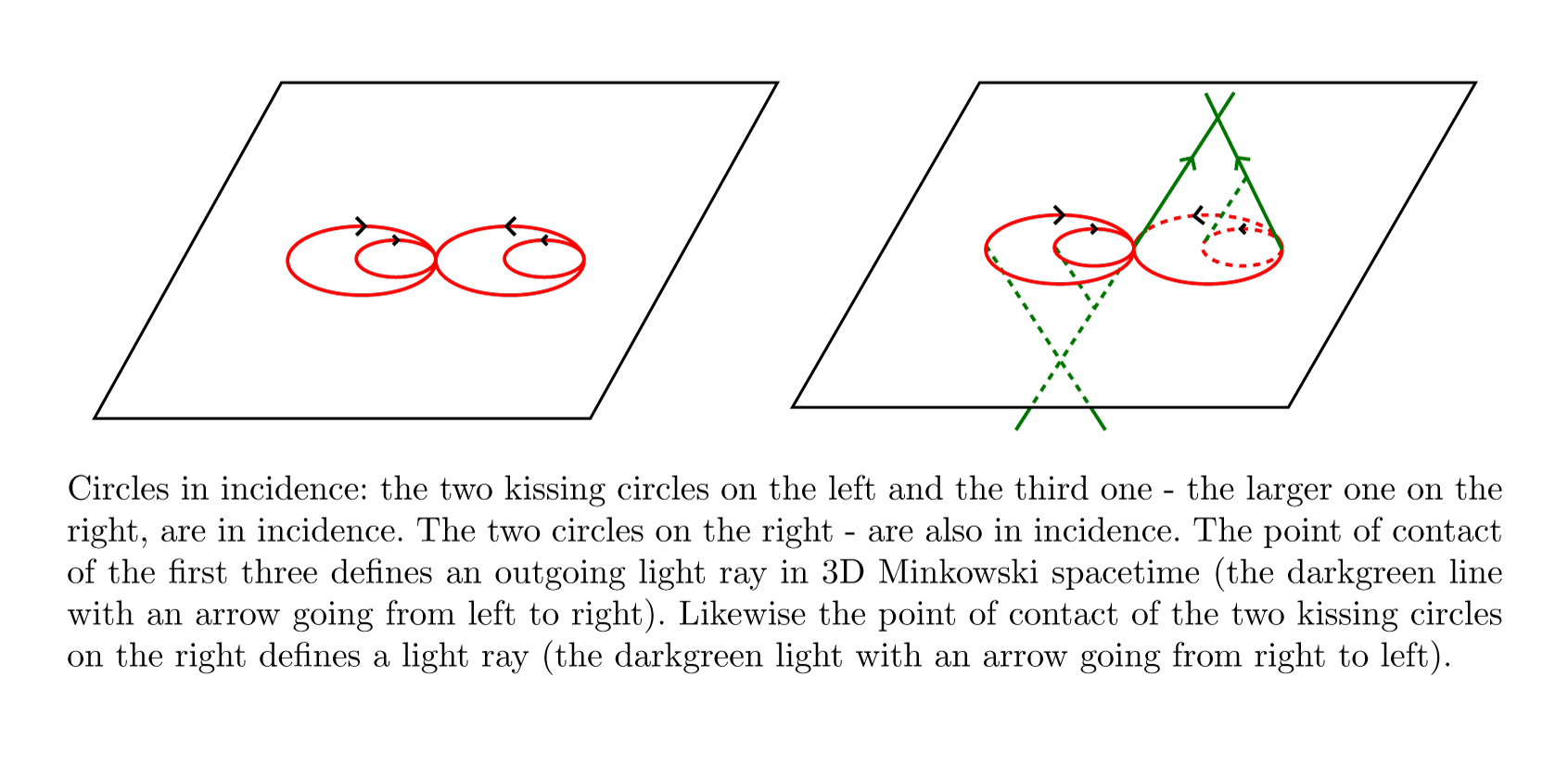}}

Indeed, parametrizing the space of circles on the plane by $(a,b,R)$, where  $(a,b)$ are the coordinates of their center in the plane, and $R$ is their (negative or positive) radius, we see that close circles corresponding to $(a,b,R)$ and $(a+\der a,b+\der b,R+\der R)$ have only one point of intersection iff the equations 
$$(x-a)^2+(y-b)^2-R^2=0\quad\&\quad (x-a-\der a)^2+(y-b-\der b)^2-(R+\der R)^2=0,$$ have a unique solution for $(x,y)$. It is only possible if and only if
$$(\der a)^2+(\der b)^2-(\der R)^2=0,$$
i.e. when the circles corresponding to $(a,b,R)$ and $(a+\der a,b+\der b,R+\der R)$ are \emph{null} separated in the Lorentzian metric $g=(\der a)^2+(\der b)^2-(\der R)^2$ on the space of all circles ${\color{red}{\bf Q}_c}$. Thus the space of all circles ${\color{red}{\bf Q}_c}$ is embedded as an open set in the projective quadric ${\color{red}{\bf Q}}$, and moreover this embedding is a \emph{conformal embedding} with a \emph{flat conformal structure} coming from the \emph{Minkowski metric} $g=(\der a)^2+(\der b)^2-(\der R)^2$.

Another, more geometric, way of seeing the conformal metric $g=(\der a)^2+(\der b)^2-(\der R)^2$ on the space ${\color{red}{\bf Q}_c}$ is to think about $(a,b)$ plane as a $R=0$ slice of $\bbR^3$ with coordinates $(a,b,R)$. This space can be uniquely equipped with the set of cones, such that each circle with center in $(a_0,b_0)$ and (positive or negative) radius $R_0$ on the $R=0$ plane is an intersection of this plane with a cone having tip at $(a_0,b_0,R_0)$. Then one declares $\bbR^3$ with such cones as a conformal 3-dimensional manifold on which these cones are light cones. By construction these cones are light cones in the metric $g=(\der a)^2+(\der b)^2-(\der R)^2$.
\subsubsection{Conformal Minkowski space in 3-dimensions is $\sog(2,3)$ symmetric} Since, following Lie, we have shown that the space ${\color{red}{\bf Q}}$ of all circles on the plane has a natural structure of 3-dimensional conformal Minkowski space which has $\sog(2,3)$ as a group of symmetries, and since ${\color{red}{\bf Q}}$ is in one to one correspondence with the base ${\color{red}Q}$ of the fibration $M\to {\color{red}Q}$, then also the space  ${\color{red}Q}$ of all integral curves of the vector field ${\color{red}X_4}$ in $M$ has $\sog(2,3)$ as a symmetry. But this is naturally associated with the configuration space $M$ of a car equipped with the geometry of an (velocity) Engel distribution with car's split. This gives an argument why the Lie algebra $\soa(2,3)$ is the algebra of infinitesimal symmetries of the car structure $(M,{\color{red}{\mathcal D}}\hspace{-0.31cm}{\color{darkgreen}{\mathcal D}}={\color{darkgreen}{\mathcal D_w}}\oplus {\color{red}{\mathcal D_g}})$.
\subsection{Geometry of 3rd order ODEs}\label{secode}
It turns out that the double fibration of the type\\
 \centerline{\includegraphics[height=5cm]{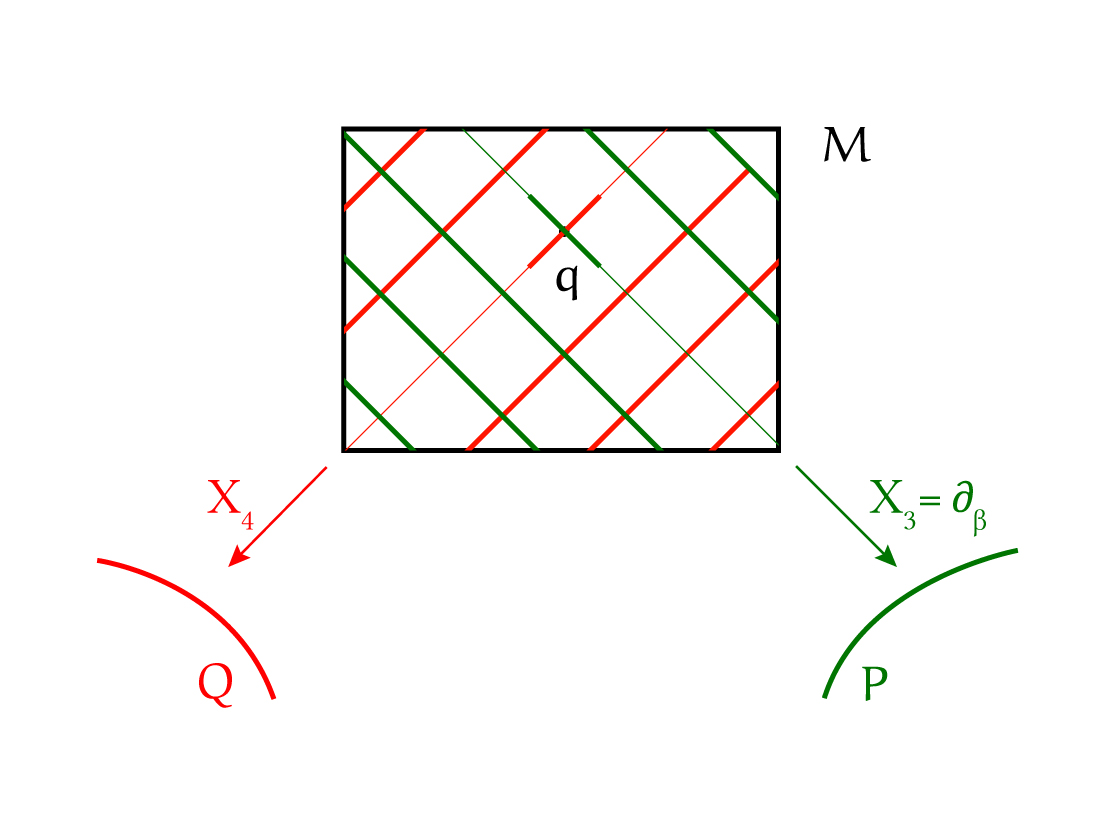}}
 is also associated with the geometry of \emph{3rd order ODEs considered modulo contact transformations of variables}. Indeed, in \cite{chern} S.S. Chern studied the geometry of an ordinary differential equation (ODE) 
 \be y'''=F(x,y,y',y'')\label{ceq}\ee considered modulo contact transformation of variables, and established that the space $M$ of second jets of the ODE, i.e. the \emph{four}-dimensional space coordinatized by the jet coordinates $(x,y,y',y'')$, is naturally equipped with two 1-dimensional foliations. These are given
 \begin{itemize}
 \item in terms of the integral curves of a vector field ${\color{darkgreen}X_3=\partial_{y''}}$ responsible for the projection $(x,y,y',y'')\to (x,y,y')$ from the space $M$ of second jets to space ${\color{darkgreen}P}$ of the first jets, and 
 \item  in terms of the \emph{total differential vector field} ${\color{red}X_4=\partial_x+y'\partial_y+y''\partial_y'+F\partial_{y''}}$ of the equation.
 \end{itemize}
 He has also shown that these two foliations on $M$ do not change when the ODE undergoes contact transformation of variables. This led him to the study of a double fibration ${\color{red}Q}\leftarrow M\rightarrow {\color{darkgreen}P}$, with the 3-dimensional space ${\color{red}Q}$ beeing the leaf space of the foliation given by ${\color{red}X_4}$.

 In this section we recall Chern's considerations, and will show their relation to the geometry of the car fibration.

The equation \eqref{ceq} can be equivalently written as a system $y'=p$, $p'=q$, $q'=F(x,y,p,q)$ and as such is defined on the space of second jets $M={\mathcal J}^2$ over the $x$-axis. This space is parameterized by $(x,y,p,q)$ and every solution to \eqref{ceq} is a curve $\gamma(t)=(x(t),y(t),p(t),q(t))$ in ${\mathcal J}^2$ such that its tangent vector $\dot{\gamma}(t)$ annihilates the contact forms
\be \omega^1=\der y-p\der x,\quad\omega^2=\der p-q\der x,\quad \omega^3=\der q-F(x,y,p,q)\der x.\label{cof1}\ee
These can be supplemented by \be \omega^4=\der x\label{cof2}\ee to a coframe on ${\mathcal J}^2$.

Chern, inspired by the earlier work of E. Cartan's \cite{car1}, (see also \cite{car2}), established that an arbitrary \emph{contact} transformation of variables of the equation \eqref{ceq} is equivalent to the following transformation of the \emph{coframe} 1-forms $(\omega^1,\omega^2,\omega^3,\omega^4)$ in ${\mathcal J}^2$:
\be\begin{aligned}
  \bma\omega^1\\\omega^2\\\omega^3\\\omega^4\ema\to \bma t_1&t_2&0&0\\t_3&t_4&0&0\\t_5&t_6&t_7&0\\t_8&t_9&0&t_{10}\ema\bma\omega^1\\\omega^2\\\omega^3\\\omega^4\ema.
  \end{aligned}
\label{trcof}\ee
Here the $t_i$ are arbitrary functions on ${\mathcal J}^2$ such that $(t_1t_4-t_2t_3)t_7t_{10}\neq 0$. Thus the local equivalence of 3-rd order ODEs, considered modulo contact transformations, got reformulated by Chern into the local equivalence of coframes \eqref{cof1}-\eqref{cof2} given modulo transformations \eqref{trcof}.

Looking at the transformation \eqref{trcof} defining a contact equivalence class of ODEs \eqref{ceq}, we see that the frame vector fields $(X_1,X_2,{\color{darkgreen}X_3},{\color{red}X_4})$, which on ${\mathcal J}^2$ are dual to $(\omega^1,\omega^2,\omega^3,\omega^4)$, $X_i\hook\omega^j=\delta_i{}^j$, are given up to the transformations
\be\begin{aligned}
  \bma X_1\\X_2\\{\color{darkgreen}X_3}\\{\color{red}X_4}\ema\to \bma *&*&*&*\\ *&*&*&*\\0&0&\tfrac{1}{t_7}&0\\ 0&0&0&\tfrac{1}{t_{10}}\ema\bma X_1\\X_2\\{\color{darkgreen}X_3}\\{\color{red}X_4}\ema.
  \end{aligned}
\label{trf}\ee
Thus a 3rd order ODE \eqref{ceq} considered modulo contact transformations distinguishes two well defined directions on ${\mathcal J}^2$. They are spanned by the respective vector fields
$${\color{darkgreen}X_3=\partial_p}\quad\mathrm{and}\quad {\color{red}X_4=\partial_x+p\partial_y+q\partial_p+F\partial_q}.$$

These in turn span a rank 2 distribution ${\color{red}{\mathcal D}}\hspace{-0.31cm}{\color{darkgreen}{\mathcal D}}=\Span_{{\mathcal F}({\mathcal J}^2)}({\color{darkgreen}X_3},{\color{red}X_4})$ which happens to be an Engel distribution. Thus we have an Engel distribution ${\color{red}{\mathcal D}}\hspace{-0.31cm}{\color{darkgreen}{\mathcal D}}$ with a natural split ${\color{red}{\mathcal D}}\hspace{-0.31cm}{\color{darkgreen}{\mathcal D}}={\color{darkgreen}{\mathcal D}_w}\oplus{\color{red}{\mathcal D}_g}$ given by ${\color{darkgreen}{\mathcal D}_w}=\Span_{{\mathcal F}({\mathcal J}^2)}({\color{darkgreen}X_3})$ and ${\color{red}{\mathcal D}_g}=\Span_{{\mathcal F}({\mathcal J}^2)}({\color{red}X_4})$. So the geometry of the jet space ${\mathcal J}^2$ with a 3rd order ODE considered modulo point transformations of variables is very much like the geometry of car's configuration space!

Can we thus associate a 3rd order ODE to the car? If so, what is the ODE?

It turns out that the car structure geometry is a special case of geometries studied by us in the paper  
\cite{bulletin}. There we considered manifolds $M$ of dimension $k+n$ and the geometry of rank $n=r+s$ distributions $\mathcal D$ on $M$ which had the split ${\mathcal D}={\mathcal D}_r\oplus{\mathcal D}_s$ onto \emph{integrable} subdistributions of respective ranks $r$ and $s$. We called such structures \emph{para-CR structures of type} $(k,r,s)$. Since rank 1 distributions are always integrable then, in this sense, the geometry of car's structure $(M,{\color{red}{\mathcal D}}\hspace{-0.31cm}{\color{darkgreen}{\mathcal D}}={\color{darkgreen}{\mathcal D}_w}\oplus{\color{red}{\mathcal D}_g})$ is a para-CR structure of type $(2,1,1)$.

Actually, in Ref. \cite{bulletin}, Sec. 4, Proposition 4.2, we have shown that the geometry of para-CR structures of type $(2,1,1)$ is the same as the geometry of 3rd order ODEs considered modulo contact transformation of variables. Thus, according to this general result, there definitely exists a contact equivalence class of 3rd order ODEs associated with a car. So what is an ODE representing this class? 

The car structure  $(M,{\color{red}{\mathcal D}}\hspace{-0.31cm}{\color{darkgreen}{\mathcal D}}={\color{darkgreen}{\mathcal D}_w}\oplus{\color{red}{\mathcal D}_g})$ defines a $G$-structure \cite{gstr} on $M$, i.e. the reduction of the structure group $\glg(4,\bbR)$ of the tangent bundle $\mathrm{T}M$ to its subgroup
$G=\{\glg(4,\bbR)\ni A: A{\color{darkgreen}X_3}=\lambda_3 {\color{darkgreen}X_3}, A{\color{red}X_4}=\lambda_4 {\color{red}X_4}\}$, preserving ${\color{red}{\mathcal D}}\hspace{-0.31cm}{\color{darkgreen}{\mathcal D}}$ and its split ${\color{red}{\mathcal D}}\hspace{-0.31cm}{\color{darkgreen}{\mathcal D}}={\color{darkgreen}{\mathcal D}_w}\oplus{\color{red}{\mathcal D}_g}$. It is more convenient to think about a $G$-structure dually: it is a $G$-subbundle of the bundle $F^*(M)$ of $\glg(4,\bbR)$-coframes of $M$. The requirement that the $G$-structure is given by the car structure $(M,{\color{red}{\mathcal D}}\hspace{-0.31cm}{\color{darkgreen}{\mathcal D}}={\color{darkgreen}{\mathcal D}_w}\oplus{\color{red}{\mathcal D}_g})$ is reflected in the $G$ transformation of coframes as follows. We first consider the coframe $(\omega^1,\omega^2,\omega^3,\omega^4)$ dual to the car frame $(X_1,X_2,{\color{darkgreen}X_3},{\color{red}X_4})$ on $M$ given in \eqref{vd2}, \eqref{vd4}. We have:
\be
\begin{aligned}
  \omega^1&=\ell^{-1}(\cos\alpha\der y-\sin\alpha\der x)\\
  \omega^2&=-\cos\beta\der\alpha-\ell^{-1}\sin\beta\Big(\cos\alpha\der x+\sin\alpha\der y\Big)\\
  \omega^3&=\der\beta\\
  \omega^4&=-\sin\beta\der\alpha+\ell^{-1}\cos\beta\Big(\cos\alpha\der x+\sin\alpha\der y\Big),
\end{aligned}
\label{cof}
\ee
and $X_i\hook\omega^j=\delta_i{~}^j$. Now, the coframe $(\omega^i)$, $i=1,2,3,4$, is given by the geometry of the car up to the transformation
\be \omega^i\to \bar{\omega}{}^i=A^i{~}_j\omega^j,\label{tra}\ee
with
\be A=(A^i{~}_j)=\bma t_1&t_2&0&0\\t_3&t_4&0&0\\t_5&t_6&t_7&0\\t_8&t_9&0&t_{10}\ema\quad\quad\mathrm{with}\quad\quad t_B\in{\mathcal F}(M),\quad\mathrm{and}\quad\mathrm{det}A\neq 0.\label{Group}\ee
The $G$-structure group $G$ of the car structure is therefore
$$G=\{A\in M_{4\times 4}(\bbR)~:~ A=\bma t_1&t_2&0&0\\t_3&t_4&0&0\\t_5&t_6&t_7&0\\t_8&t_9&0&t_{10}\ema\,\,\mathrm{with}\,\, t_B\in\bbR,\,\,\mathrm{and}\,\,\mathrm{det}A\neq 0\}.$$
  We now use transformations \eqref{tra}-\eqref{Group} to bring the coframe forms \eqref{cof} to a form which is convenient to see a 3rd order ODE related to the car's geometry.

Taking \be A_1=\bma \ell\sec\alpha&0&0&0\\0&1&0&0\\0&0&1&0\\0&0&0&1\ema\label{a1}\ee
we bring $\omega^1$ into the form \be\omega^1=\der y-\tan\alpha\der x.\label{om1}\ee
Now we observe that $$\omega^2=-\cos\beta\cos^2\alpha\Big(\der\tan\alpha+\ell^{-1}\tan\beta\sec^3\alpha\der x\Big)-\ell^{-1}\sin\beta\sin\alpha\omega^1,$$
where we have used the new $\omega^1$ given by \eqref{om1}. This means that by taking
\be A_2=\bma 1&0&0&0\\-\ell^{-1}\tan\beta\tan\alpha\sec\alpha&-\sec\beta\sec^2\alpha&0&0\\0&0&1&0\\0&0&0&1\ema\label{a2}\ee
we can bring the coframe 1-form $\omega^2$ into the form
\be\omega^2=\der\tan\alpha+\ell^{-1}\tan\beta\sec^3\alpha\der x.\label{om2}\ee
We further observe that
$$\omega^3=-\ell^{-1}\sec^3\alpha\sec^2\beta\Big(-\der\big(\ell^{-1}\sec^3\alpha\tan\beta\big)-3\ell^{-2}\sec^5\alpha\sin\alpha\tan^2\beta\der x\Big)-\tfrac34\sin2\alpha\sin2\beta\omega^2,$$
where we have used the new $\omega^2$ given by \eqref{om2}. This means that by means of the matrix
\be A_3=\bma 1&0&0&0\\0&1&0&0\\0&-3\ell^{-1}\sec\alpha\tan\alpha\tan\beta&-\ell^{-1}\sec^3\alpha\sec^2\beta&0\\0&0&0&1\ema\label{a3}\ee
we can bring the 1-form $\omega^3$ into the form
\be\omega^3=-\der\big(\ell^{-1}\sec^3\alpha\tan\beta\big)-3\ell^{-2}\sec^5\alpha\sin\alpha\tan^2\beta\der x.\label{om3}\ee
Finally, we also see that
$$\omega^4=\ell^{-1}\sec\alpha\sec\beta\der x-\cos^2\alpha\sin\beta\omega^2+\ell^{-1}\cos\beta\sin\alpha\omega^1,$$
with $\omega^1$ and $\omega^2$ as in \eqref{om1}, \eqref{om2}, which shows that the matrix
\be A_4=\bma 1&0&0&0\\0&1&0&0\\0&0&1&0\\-\tfrac12\cos^2\beta\sin2\alpha&\tfrac12\ell\cos^3\alpha\sin2\beta&0&\ell\cos\alpha\cos\beta\ema\label{a4}\ee
brings the form $\omega^4$ into
\be\omega^4=\der x.\label{om4}\ee
Summarizing what we have obtained so far we note that by a linear transformation
$$A=A_4A_3A_2A_1,$$
with $A_i$ as in \eqref{a1}, \eqref{a2}, \eqref{a3}, \eqref{a4}, which is of the form of \eqref{Group}, we can bring the car coframe \eqref{cof} to the $G$-equivalent coframe
\be\begin{aligned}
  \omega^1&=\der y-\tan\alpha\der x\\
  \omega^2&=\der\tan\alpha+\ell^{-1}\tan\beta\sec^3\alpha\der x\\
  \omega^3&=-\der\big(\ell^{-1}\sec^3\alpha\tan\beta\big)-3\ell^{-2}\sec^5\alpha\sin\alpha\tan^2\beta\der x\\
  \omega^4&=\der x.
\end{aligned}
\label{cofe}\ee
Now we introduce the new coordinates $(x,y,p,q)$ on $M$ related to the coordinates $(x,y,\alpha,\beta)$ via
$$p=\tan\alpha,\quad\quad q=-\ell^{-1}\tan\beta\sec^3\alpha.$$
In these new coordinates the coframe 1-forms \eqref{cofe} read:
$$\begin{aligned}
  \omega^1&=\der y-p\der x\\
  \omega^2&=\der p-q\der x\\
  \omega^3&=\der q-F(x,y,p,q)\der x\\
  \omega^4&=\der x,
\end{aligned}
$$
with
$$F=3\ell^{-2}\sec^5\alpha\sin\alpha\tan^2\beta=\frac{3pq^2}{1+p^2}.$$
Thus the car structure can equivalently be described in terms of coordinates $(x,y,p,q)$ with the adapted coframe 1-forms
\be\begin{aligned}
  \omega^1&=\der y-p\der x\\
  \omega^2&=\der p-q\der x\\
  \omega^3&=\der q-\frac{3pq^2}{1+p^2}\der x\\
  \omega^4&=\der x.
\end{aligned}
\label{coffe}\ee
The car velocity distribution $${\color{red}{\mathcal D}}\hspace{-0.31cm}{\color{darkgreen}{\mathcal D}}=\Span_{{\mathcal F}(M)}({\color{darkgreen}X_3},{\color{red}X_4})$$ is in these coordinates spanned by the vector fields
\be
{\color{darkgreen}X_3=\partial_q}\quad\mathrm{and}\quad {\color{red}X_4=\partial_x+p\partial_y+q\partial_p+\frac{3pq^2}{1+p^2}\partial_q}.\label{newvec}\ee
They form a part of a frame $(X_1,X_2,{\color{darkgreen}X_3},{\color{red}X_4})$ dual to $(\omega^1,\omega^2,\omega^3,\omega^4)$ given by \eqref{coffe}. The
`steering wheel'-`gas' split, $${\color{red}{\mathcal D}}\hspace{-0.31cm}{\color{darkgreen}{\mathcal D}}={\color{darkgreen}{\mathcal D}_w}\oplus{\color{red}{\mathcal D}_g},$$ is given by $${\color{darkgreen}{\mathcal D}_w}=\Span_{{\mathcal F}(M)}({\color{darkgreen}X_3})\quad\mathrm{and}\quad{\color{red}{\mathcal D}_g}=\Span_{{\mathcal F}(M)}({\color{red}X_4}).$$
Since the coframe 1-forms \eqref{coffe} are just the standard \emph{contact forms on the bundle of second jets} ${\mathcal J}^2$ with the standard jet coordinates $(x,y,p=y',q=y'')$ as in \eqref{cof1}-\eqref{cof2}, we recognize here \emph{the third ODE} $y'''=F(x,y,y',y'')$, with $F=\frac{3pq^2}{1+p^2}$. The possible transformations \eqref{tra}-\eqref{Group} of these forms, are equivalent to the \emph{contact transformations} of variables for this equation (see \cite{bulletin}, Sec. 4). Thus, the geometry of the car structure $(M,{\color{red}{\mathcal D}}\hspace{-0.31cm}{\color{darkgreen}{\mathcal D}}={\color{darkgreen}{\mathcal D}_w}\oplus{\color{red}{\mathcal D}_g})$ is locally diffeomorphically equivalent to the local differential geometry of the 3rd order ODE
\be y'''=\frac{3y'y''{}^2}{1+y'{}^2}\label{code}\ee
considered modulo contact transformation of variables.

\emph{What is this equation?}
This is the equation whose graphs of \emph{general} solutions $(x,y(x))$ describe all circles on the plane $(x,y)$. Indeed one can easilly check that the general solution to \eqref{code} is given by
  $$\nu(x^2+y^2)-2\xi x-2\eta y+\mu=0,$$
where $\nu,\xi,\eta,\mu$ are real constants. Since this formula is projective, the space of solutions ${\color{red}{\bf Q}}$ is 3-dimensional.
Taking $\nu=1$ we get the space ${\color{red}{\bf Q}_c}$ of all circles on the plane (with radius $R=\sqrt{\eta^2+\xi^2-\mu}$, centered at $x=a=\xi$ and $y=b=\eta$), taking $\nu=0$ we get the space ${\color{red}{\bf Q}_\ell}$ of all lines in the plane, and taking $\nu=1$ and $\eta^2+\xi^2=\mu$ we get the space ${\color{red}{\bf Q}_p}$ of all points in the plane. 

\emph{What is the relation of the circles $(x^2+y^2)-2\xi x-2\eta y+\mu=0$ and the lines $2\xi x+2\eta y-\mu=0$ to the car movement?}  
By construction, the vector field ${\color{red}X_4}$ in \eqref{newvec} differs from the vector field ${\color{red}X_4}$ in \eqref{vd2} by rescaling. Thus, \emph{modulo a reparametrization}, both of them have \emph{the same} integral curves in $M$. We know that the curves defined by ${\color{red}X_4}$ from \eqref{vd2} are \emph{helixes} ($\beta_0\neq 0$) or \emph{straight lines} ($\beta_0=0$), which when projected on the $(x,y)$ plane, are \emph{circles} or \emph{straight lines} there. Likewise the integral curves of ${\color{red}X_4}$ from \eqref{newvec} are helixes or straight lines which project to the circles or straight lines in the $(x,y)$ plane. To see this one considers a curve ${\color{red}q_4(t)=(x(t),y(t),p(t),q(t))}$ in $M$ such that ${\color{red}\dot{q}_4}$ is tangent to ${\color{red}X_4}$ from \eqref{newvec}. It satisfies the system of ODEs $(\dot{x},\dot{y},\dot{p},\dot{q})=(1,p,q,\frac{3pq^2}{1+p^2})$. This means that $x=t$, $p=\dot{y}$, $q=\dot{p}=\ddot{y}$, and finally $\dot{q}=\dddot{y}=\frac{3\dot{y}\ddot{y}^2}{1+\dot{y}^2}.$ Thus, the graphs of solutions $y=y(t)$ of the last equation in the plane $(x=t,y)$, which are \emph{circles} or \emph{straight lines}, are just the circles or straight lines which the rear wheels of the car are performing in the plane $(x,y)$ when the driver of a car applies a primitive `gas control' only.

\subsection{Contact projective geometry on ${\color{darkgreen}P}$}
We now pass to analyse the geometry of ${\color{darkgreen}P}$, i.e. the base of the fibration $M\to{\color{darkgreen}P}$, whose fibers are the steering wheel trajectories  generated by the steering wheel vector field ${\color{darkgreen}X_3}$ on the car's configuration space $M$. So what is the geometry on ${\color{darkgreen}P}$?

To answer this question let us start with the interpretation of the configuration space $M$ of the car as the second jet space for the car's ODE $y'''=\frac{3y'y''{}^2}{1+y'{}^2}$. In this interpretation, the unparametrized integral curves of ${\color{darkgreen}X_3=\partial_\beta}$ are the same as the unparametrized integral curves of ${\color{darkgreen}X_3=\partial_q}$, and they constitute natural fibres of the fibration $\pi:M={\mathcal J}^2\to {\color{darkgreen}P={\mathcal J}^1}$ of the second jet space ${\mathcal J}^2$ with coordinates $(x,y,p,q)$ over the first jet space $\color{darkgreen}{\mathcal J}^1$ with coordinates $(x,y,p)$. Consider now the trajectories of ${\color{red}X_4}$, which in the second jet interpretation of $M$, are just curves ${\color{red}(x,y,p,q)=(x,y(x),y'(x),y''(x))}$ in ${\mathcal J}^2$ corresponding to \emph{solutions} $y=y(x)$ of the ODE \eqref{code}. There is a natural projection $\pi({\color{red}(x,y(x),y'(x),y''(x))})=\color{red}(x,y(x),y'(x))$ of these curves to the 3-dimensional space $\color{darkgreen}{\mathcal J}^1$ of the first jets. The important observation is that these projected curves ${\color{red}(x,y,p)=(x,y(x),y'(x))}$ in $\color{darkgreen}{\mathcal J}^1$, as curves corresponding to the solutions of \eqref{code}, are \emph{always tangent to the contact distribution} ${\mathcal C}=\{X\in \Gamma(\mathrm{T}{\mathcal J}^1)~:~ X\hook(\der y-p\der x)=0\}$, which is a \emph{natural} structure on ${\mathcal J}^1$. Moreover, since we have a solution to \eqref{code} for every choice of initial conditions $y(x_0)=y_0$, $y'(x_0)=p_0$, then at every point $(x_0,y_0,p_0)$ in $\color{darkgreen}{\mathcal J}^1$ the projection  $\pi((x,y(x),y'(x),y''(x)))$ defines a curve tangent to $\mathcal C$ in \emph{every} direction of $\mathcal C$. It follows that the projections  $\pi({\color{red}(x,y(x),y'(x),y''(x))})$ of solution curves from ${\mathcal J}^2$ to ${\mathcal J}^1$ can be considered as \emph{geodesics of} a certain class of \emph{torsion free connections} on ${\color{darkgreen}P={\mathcal J}^1}$.

Indeed, consider a curve $\gamma(t)=x(t)\partial_x+y(t)\partial_y+p(t)\partial_p$ tangent to a distribution $\mathcal C$ in $\color{darkgreen}{\mathcal J}^1$, and a frame $(Z_1,Z_2,Z_3)$ in $\color{darkgreen}{\mathcal J}^1$ with
$$Z_1=\partial_y,\quad Z_2=\partial_x+\partial_y,\quad Z_3=\partial_p.$$
Since ${\mathcal C}=\Span_{{\mathcal F}({\mathcal J}^1)}(Z_2,Z_3)$, the velocity of this curve, $$\dot{\gamma}=\dot{x}\partial_x+\dot{y}\partial_y+\dot{p}\partial_p=\dot{\gamma}{}^1Z_1+\dot{\gamma}{}^2Z_2+\dot{\gamma}{}^3Z_3=(\dot{y}-p\dot{x})\partial_y+\dot{x}Z_2+\dot{p}Z_3,$$ has the following components in the frame $(Z_1,Z_2,Z_3)$:
\be \dot{\gamma}{}^1=\dot{y}-p\dot{x}=0, \quad \dot{\gamma}{}^2=\dot{x},\quad \dot{\gamma}{}^3=\dot{p}.\label{veldot}\ee
If the curve $\gamma(t)$ is a geodesic of a torsion free connection, there should exist functional coefficients $\Gamma^i{}_{jk}=\Gamma^i{}_{kj}$ - the connection coefficients in the frame $(Z_1,Z_2,Z_3)$ - such that
$$\dot{\gamma}{}^i+\Gamma^i{}_{jk}\gamma^j\gamma^k=0.$$
Thus, to interprete  $\gamma(t)$ as a geodesic it is enough to find $\Gamma^i{}_{jk}=\Gamma^i{}_{kj}$ such that
\be \ddot{x}+\Gamma^2{}_{22}\dot{x}{}^2+2\Gamma^2{}_{23}\dot{x}\dot{p}+\Gamma^2{}_{33}\dot{p}{}^2=0\quad\quad\&\quad\quad\ddot{p}+\Gamma^3{}_{22}\dot{x}{}^2+2\Gamma^3{}_{23}\dot{x}\dot{p}+\Gamma^3{}_{33}\dot{p}{}^2=0.\label{geoo}\ee
For this we eliminate $t$ from both of these equations, by parametrizing $y=y(t)$ and $p=p(t)$ by $x$. Because of the first equation in \eqref{veldot} we have
$$p=\frac{\dot{y}}{\dot{x}}=\frac{\der y}{\der x}=y',\quad \dot{p}=\dot{x}y'',\quad\ddot{p}=\ddot{x}y''+\dot{x}{}^2y''',$$
and the last two of these equations compared with the second equation in \eqref{geoo} shows that
$$-(\Gamma^3{}_{22}\dot{x}{}^2+2\Gamma^3{}_{23}\dot{x}{}^2y''+\Gamma^3{}_{33}\dot{x}{}^2y''{}^2)=-y''(\Gamma^2{}_{22}\dot{x}{}^2+2\Gamma^2{}_{23}\dot{x}{}^2y''+\Gamma^2{}_{33}\dot{x}{}^2y''{}^2)+\dot{x}{}^2y'''.$$
  Simplifying, we get:
  $$y'''=\Gamma^2{}_{33}y''{}^3+(2\Gamma^2_{23}-\Gamma^3{}_{33})y''{}^2+(\Gamma^2{}_{22}-2\Gamma^3{}_{23})y''-\Gamma^3{}_{22},$$
  where $\Gamma^i{}_{jk}$ are functions of $x,y$ and $p=y'(x)$ only.

  Thus, for an equation $y'''=F(x,y,y',y'')$ to define on the space of first jets ${\mathcal J}^1$ a structure of a contact manifold with geodesics passing through every point in every direction and such that they are tangent to the contact distribution, it is neccessary that the function $F=F(x,y,y',y'')$ is a polynomial of at most 3rd order in the variable $y''$. It follows that this condition for $F$ is also sufficient for getting such a structure on ${\mathcal J}^1$.

  Since the car's structure equation $y'''=\frac{3y'y''{}^2}{1+y'{}^2}$ depends on $y''$ quadratically, this implies that its 3-dimensional space ${\color{darkgreen}P}$ i.e. its space of first jets $\color{darkgreen}{\mathcal J}^1$ is naturally equipped with the structure as in the following definition \cite{fox}.

  \begin{definition}
A \emph{contact projective structure} on the first jet space ${\mathcal J}^1$ is given
by the following data.
\begin{itemize}
\item The contact distribution $\mathcal C$, that is the distribution annihilated by
$\omega^1= dy-pdx$.
\item A family of unparameterized curves everywhere tangent to ${\mathcal C}$ and such that:
  \begin{itemize}
    \item for a given point and a direction in $\mathcal C$ there is exactly one curve passing
      through that point and tangent to that direction,
      \item curves of the family
        are among unparameterized geodesics for some linear connection on ${\mathcal J}^1$.
  \end{itemize}
  \end{itemize}
\end{definition}
To make the statement above the definition more explicit, we argue as follows:

We have the fibration $M\to{\color{darkgreen}P}$, which on the one hand is a fibration of the second jet space $M={\mathcal J}^2$ of a contact equivalence class of ODEs  $y'''=\frac{3y'y"{}^2}{1+y'{}^2}$ over the space ${\color{darkgreen}P}={\mathcal J}^1$, and on the other hand the car fibration $M\to{\color{darkgreen}P}$ of the configuration space $M$ of a car and the space ${\color{darkgreen}P}$ of the possible movements of the car modulo the moves of a steering wheel. As we explained in Section \ref{secode} there is a natural bundle isomorphism between the car configuration space and the space of second jets of the contact equivalence classes of ODEs represented by $y'''=\frac{3y'y"{}^2}{1+y'{}^2}$, making an equivalence between the car's Engel geometry with a split and the contact geometry of this ODE.  
Since the ODE $y'''=\frac{3y'y"{}^2}{1+y'{}^2}$ has only quadratic dependence on $y''$ it belongs to the class of ODEs $y'''=A_3y''{}^3+A_2y''{}^2+A_1y''+A_0$. Thus the car's ODE first jet space ${\mathcal J}^1$ has a natural \emph{contact projective structure}. This, via the bundle isomorphism ${\mathcal J}^2\to{\color{darkgreen}P}$, induced by the isomorphism between geometries on the car's configuration space and the bundle of the second jets, shows that the \emph{car's 3-dimensional space} ${\color{darkgreen}P}$ of leaves generated by ${\color{darkgreen}X_3=\partial_\beta}$ has a natural \emph{contact projective structure}. Such structures were in particular studied in \cite{fox,godphd,gn}. Since the car's ODE $y'''=\frac{3y'y"{}^2}{1+y'{}^2}$  is contact equivalent to $y'''=0$, it follows from these studies that this contact projective structure is \emph{Cartan flat}. More precisely we have the following theorem:
\begin{theorem} The car's Engel structure with a split $(M,{\color{red}{\mathcal D}}\hspace{-0.31cm}{\color{darkgreen}{\mathcal D}}={\color{darkgreen}{\mathcal D}_w}\oplus{\color{red}{\mathcal D}_g})$ induces a natural contact projective structure on the car's space ${\color{darkgreen}P}$ of all possible positions of a car considered modulo orientation of the front wheels. This contact projective structure has a 10-dimensional Lie algebra of symmetries, which is isomorphic to the simple Lie algebra $\spa(2,\bbR)$. It is flat in the sense of having vanishing curvature of the natural normal $\spa(2,\bbR)$-valued Cartan connection uniquely defined by this contact projective structure.
  \end{theorem}

Since $\spa(2,\bbR)$ is isomorphic to $\soa(2,3)$ (see Section \ref{sec511}) we again have an indication why the geometry of car's configuration space $M$ has $\soa(2,3)$ as its local symmetry.
\subsection{Chern's double fibration ${\color{red}Q}\leftarrow M\to{\color{darkgreen}P}$, the geometries on ${\color{red}Q}$ and ${\color{darkgreen}P}$ and a problem about a car on a curved terrain}
If somebody inspired by this article would like to \emph{curve} the geometry of a car, she will find usefull the following information about the geometry of general third order ODEs, $y'''=F(x,y,y',y'')$, considered modulo contact transformations.\\
 \centerline{\includegraphics[height=6cm]{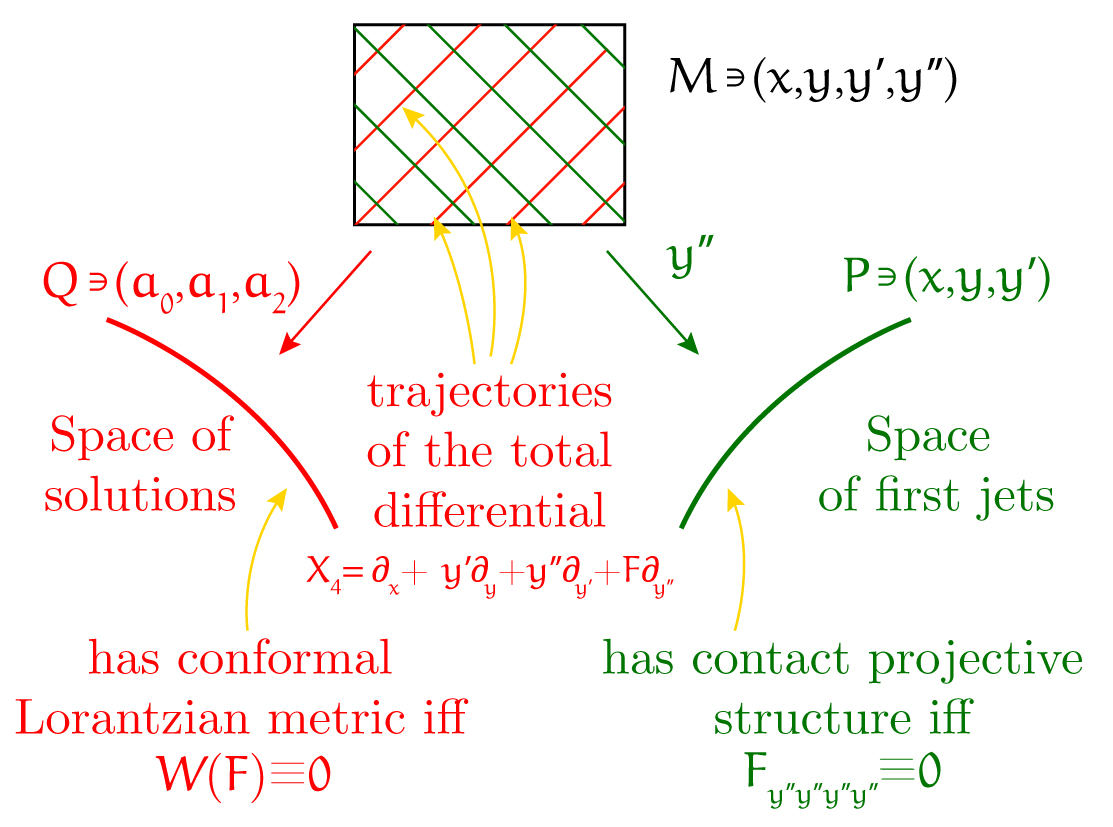}}
 As we mentioned in Section \ref{secode} Chern in 1940 noticed the above double fibration ${\color{red}Q}\leftarrow M\to \color{darkgreen}P$ for any contact equivalence class of third order ODEs. If the class is defined by the equation $y'''=F(x,y,y',y'')$, and if the general solution of the defining the equation is written as $y=y(x,a_1,a_2,a_3)$, where $a_1$, $a_2$, $a_3$ are the three constants of integration, then the base space ${\color{red}Q}$ is the leaf space of the total differential ${\color{red}X_4=\partial_x+y'\partial_y+y''\partial_{y'}+F\partial_{y''}}$, which is parameterized by $(a_1,a_2,a_3)$, and the base space ${\color{red}P}$ is the space of first jets - the leaf space of the integral curves of the vector field ${\color{darkgreen}X_3=\partial_{y''}}$, which is parameterized by $(x,y,y')$. We emphasize that this double fibration exists for any choice of the function $F$, and in turn is associated with \emph{any} contact equivalence class of 3rd order ODEs. However, and this is the main observation of S.S. Chern in \cite{chern}, the space of solutions ${\color{red}Q}$ has a natural conformal Lorentzian geometry on it, and/or the first jet space ${\color{darkgreen}P}$ has a natural contact projective structure on it, if and only if the function $F$ satisfies certain conditions, which are invariant with respect to contact change of the variables of the equation.

 We have the following theorem \cite{chern,godphd,gn}.
 \begin{theorem}
   The space ${\color{red}Q}$ in Chern's double fibration ${\color{red}Q}\leftarrow M\to \color{darkgreen}P$ associated with a contact equivalence class of ODEs $y'''=F(x,y,y',y'')$ has a natural conformal Lorentzian structure on it, if and only if the W\"unschmann invariant
$$W[F]=
9\,{\color{red}X_4}({\color{red}X_4}(\,{\color{darkgreen}X_3}(\,F\,))
-
27\,{\color{red}X_4}(\,F_{y'}\,)
-
18\,{\color{darkgreen}X_3}(\,F\,)\,{\color{red}X_4}(\,{\color{darkgreen}X_3}(\, F\,)\,)
+
18\,{\color{darkgreen}X_3}(\,F\,)\,F_{y'}
+
4\,{\color{darkgreen}X_3}(\,F\,)^3
+
54\,F_y
$$
identically vanishes for $F$.

Similarly,  the space ${\color{darkgreen}P}$ in the Chern's double fibration ${\color{red}Q}\leftarrow M\to \color{darkgreen}P$ associated with a contact equivalence class of ODEs $y'''=F(x,y,y',y'')$ has a natural contact projective structure on it, if and only if the Chern invariant
$$C[F]={\color{darkgreen}X_3}(\,{\color{darkgreen}X_3}(\,{\color{darkgreen}X_3}(\,{\color{darkgreen}X_3}(\,F\,)\,)\,)\,)$$
identically vanishes for $F$.

Here, ${\color{darkgreen}X_3=\partial_{y''}}$, ${\color{red}X_4=\partial_x+y'\partial_y+y''\partial_{y'}+F\partial_{y''}}$, $F_y=\frac{\partial F}{\partial y}$ and $F_{y'}=\frac{\partial F}{\partial y'}$.
   \end{theorem}

 In the car's fibration we have $F=\frac{3y'y"{}^2}{1+y'{}^2}$. This function has $W[F]\equiv C[F]\equiv 0$. Thus the car fibration has a (flat) conformal structure on ${\color{red}Q}$ and a (flat) contact projective structure on ${\color{darkgreen}P}$. This provoks the following (open)  problem.

\noindent
{\bf Problem.} \emph{Generalize the car setting enabling the car to move on a curved terrain. This should lead to a nonflat Engel structure with a split $(M,{\color{red}{\mathcal D}}\hspace{-0.31cm}{\color{darkgreen}{\mathcal D}}={\color{darkgreen}{\mathcal D}_w}\oplus{\color{red}{\mathcal D}_g})$ on the car's configuration space. Characterize, in terms of Chern's invariants $W[F]$, $C[F]$, and possibly their derivatives, those Engel structures with a split, which are configuration space structures of cars on curved terrains. Which of the two geometries: the conformal Lorentzian one, or the contact projective one will survive for a car on a general terrain? Perhaps none?}
 
\section{Lie's correspondence}
\subsection{Lagrangian planes in $\bbR^4$ and oriented circles in the plane}\label{lsg}
It was S. Lie who understood the geometry of the projective quadric $\color{red}{\bf Q}$, as in \eqref{abcd}, in terms of the geometry of Lagrangian planes in a real 4-dimensional vector space. (see \cite{bryant,helgason} for more details). To talk about \emph{Lagrangian planes} we need to have a \emph{real 4-dimensional vector space} $V$ and a \emph{symplectic form} in $V$, i.e. a 2-form $\omega\in\bigwedge^2V^*$ such that $\omega\wedge\omega\neq 0.$ Now, a 2-plane $q=\Span(Y_1,Y_2)$, with $Y_1,Y_2\in V$ and $Y_1\dz Y_2\neq 0$, is Lagrangian in $V$ if and only if $\omega(Y_1,Y_2)=0$.

Given a symplectic form $\omega$ in $V$ we consider the 5-dimensional vector space $\omega^\perp \subset\bigwedge^2V$ consisting of elements $Y\in \bigwedge^2V$ annihilating $\omega$:
$$\omega^\perp=\{\textstyle{\bigwedge^2}V\ni Y~:~Y\hook\omega=0\}.$$
It is now convenient to introduce a basis $(e_1,e_2,e_3,e_4)$ in $V$, such that the symplectic form $\omega$ reads as
\be \omega=e^1\dz e^4+e^2\dz e^3,\label{sympf}\ee
in its dual cobasis $(e^1,e^2,e^3,e^4)$, $e_i\hook e^j=\delta_i{}^j$, in $V^*$.
Then the most general element $Y\in\omega^\perp$ is:
\be Y=(\eta+\zeta)\,e_1\dz e_2+\mu\, e_1\dz e_3+\nu\, e_4\dz e_2+(\eta-\zeta)\,e_4\dz e_3+\xi\,(e_1\dz e_4-e_2\dz e_3),\label{iX}\ee
where $(\xi,\eta,\zeta,\mu,\nu)\in \bbR^5$.

We now ask the question as to when such $Y$ is a \emph{simple} bivector. Recall that an element $0\neq Y$ of $\bigwedge^2V$  is \emph{simple} if and only if $Y\dz Y=0$. In such case there exist vectors $Y_1$ and $Y_2$ in $V$ such that $Y=Y_1\dz Y_2$. Thus such $Y$ defines a 2-plane
$$q=\Span_\bbR(Y_1,Y_2),$$
in $V$.
If in addition, a \emph{simple} $Y$ belongs to the 5-dimensional subspace $\omega^\perp$, then its \emph{direction},
$${\color{red}\dr(Y)}:=\{\lambda Y,~\lambda\in\bbR\},$$
defines a 2-plane which is \emph{Lagrangian}. It turns out that \emph{every} Lagrangian 2-plane in $V$ is defined in terms of $0\neq Y\in\omega^\perp$ such that $Y\dz Y=0$.  

Simple algebra applied to a generic $Y\in\omega^\perp$ as in \eqref{iX}  gives:
\be Y\dz Y=2(\zeta^2-\eta^2+\mu\nu-\xi^2)e_1\dz e_2\dz e_3\dz e_4=-\tfrac12 Q(\xi,\eta,\zeta,\mu,\nu)e_1\dz e_2\dz e_3\dz e_4.\label{qfa}\ee
Note the appearence of the quadratic form \eqref{qform} in this formula! Thus such an $Y$ is simple, $Y\dz Y=0$, if and only if the quintuple $[\xi:\eta:\zeta:\mu:\nu]$ belongs to the projective quadric ${\color{red}{\bf Q}}$ considered in Section \ref{sec41}.
 Now, let us define
 $${\color{red}{\bf Q}'}=\{P(\textstyle{\bigwedge^2}V)\ni {\color{red}\dr(Y)}~:~ Y\hook\omega=0\quad\&\quad Y\dz Y=0\},$$
 where, as it is customary, we denoted the \emph{projectivization} of $\bigwedge^2V$ by $P(\bigwedge^2V)$.

 Since $Y\dz Y=0$ for $Y\in\omega^\perp$ is equivalent to $\zeta^2-\eta^2+\mu\nu-\xi^2=0$ for $[\xi:\eta:\zeta:\mu:\nu]\in\bbR P^4$, then  
$${\color{red}{\bf Q}'}=\{{\color{red}\dr(Y)}~:~ Y~\mathrm{as~in~\eqref{iX}}~\mathrm{with}~ [\xi:\eta:\zeta:\mu:\nu]\in {\color{red}{\bf Q}}\}.$$
This in turn establishes a diffeomorphism between ${\color{red}{\bf Q}}$ and the space of all Lagrangian 2-planes in $V$. With some abuse of notation we will denote this space also by ${\color{red}{\bf Q}'}$,
$${\color{red}{\bf Q}'}=\{\mathrm{set~of~all~Lagrangian~2}-\mathrm{planes~in}~(V,\omega)\}.$$

Let us now parametrize those ${\color{red}\dr(Y)}$ in ${\color{red}{\bf Q}'}$ that correspond to all circles with a \emph{finite} radius in the plane. Since such circles are points of the set ${\color{red}{\bf Q}_c}\subset{\color{red}{\bf Q}}$, with $\nu\neq 0$, we can conveniently parametrize them by $\nu=1$, $\mu=\xi^2+\eta^2-\zeta^2$. Thus, the corresponding bivectors ${\color{red}\dr(Y)}$ in ${\color{red}{\bf Q}'}$ may be represented by
\be Y=(\eta+\zeta)\,e_1\dz e_2+(\xi^2+\eta^2-\zeta^2)\, e_1\dz e_3+ e_4\dz e_2+(\eta-\zeta)\,e_4\dz e_3+\xi\,(e_1\dz e_4-e_2\dz e_3),\label{qc}\ee
or what is the same by $Y=\Big(\,(\eta+\zeta)\,e_1+e_4+\xi\,e_3\,\Big)\dz\Big(-\xi\,e_1+e_2+(\eta-\zeta)\,e_3\,\Big)$.
Thus in the 3-dimensional space ${\color{red}{\bf Q}'}$ there is an open set ${\color{red}{\bf Q}'_c}$ of bivectors $Y$ given by \eqref{qc}. This set, in turn, is diffeomorphic to the space of all Lagrangian 2-planes
\be q(\xi,\eta,\zeta)=\Span_\bbR\big(Y_1,Y_2\big),\label{qc1}\ee
spanned by \be Y_1=(\eta+\zeta)\,e_1+e_4+\xi\,e_3\quad \&\quad Y_2=-\xi\,e_1+e_2+(\eta-\zeta)\,e_3.\label{qc2}\ee
Again, with some abuse, we denote this space by ${\color{red}{\bf Q}'_c}$.

In ${\color{red}{\bf Q}_c}$ we had a nice interpretation of the incidence between two points (circles): two close circles in ${\color{red}{\bf Q}_c}$ were incident if they were tangent to each other. The natural \emph{incidence between the points of} ${\color{red}{\bf Q}'_c}$, i.e. between two close Lagrangian 2-planes in $V$, \emph{is their intersection along a line}. Let us see what such an incidence means:

If we take a Lagrangian 2-plane $q(\xi,\eta,\zeta)$ and its close neighbour
$q(\xi+\der\xi,\eta+\der\eta,\zeta+\der\zeta)$, then they intersect in a line iff their corresponding bivectors
$$Y=\Big(\,(\eta+\zeta)\,e_1+e_4+\xi\,e_3\,\Big)\dz\Big(-\xi\,e_1+\,e_2+(\eta-\zeta)\,e_3\,\Big)$$ and $$Y+\der Y=\Big(\,(\eta+\zeta+\der\eta+\der\zeta)\,e_1+e_4+(\xi+\der\xi)\,e_3\,\Big)\dz\Big(-(\xi+\der\xi)\,e_1+\,e_2+(\eta-\zeta+\der\eta-\der\zeta)\,e_3\,\Big)$$ satisfy
$$Y\dz(Y+\der Y)=0.$$
A short algebra shows that
$$Y\dz(Y+\der Y)=\big((\der\eta)^2+(\der\xi)^2-(\der\zeta)^2\big)e_1\dz e_2\dz e_3\dz e_4.$$
Hence the two Lagrangian planes from ${\color{red}{\bf Q}'_c}$ intersect in a line if and only if the connecting vector $(\der\xi,\der\eta,\der\zeta)$ between the points $(\xi,\eta,\zeta)$ and $(\xi+\der\xi,\eta+\der\eta,\zeta+\der\zeta)$ in ${\color{red}{\bf Q}'_c}$ is null in the 3-dimensional Minkowski metric
$g=(\der\eta)^2+(\der\xi)^2-(\der\zeta)^2$. Comparing with \eqref{abc} we see that in the present parametrization of ${\color{red}{\bf Q}_c}$, we have
$$\xi=a,\quad \eta=b \quad\mathrm{and}\quad \zeta=R.$$
Hence $g=(\der a)^2+(\der b)^2-(\der R)^2$, and the condition that two neighbouring Lagrangian planes from ${\color{red}{\bf Q}'_c}$ intersect in a line in $V$ is then equivalent to the condition that the corresponding neighbouring circles from ${\color{red}{\bf Q}_c}$ are kissing each other in the plane $(x,y)$. This is the essence of \emph{Lie's observation }:

\emph{Tangent circles in $\bbR^2$ with orientations as in moving gears correspond to Lagrangian planes in $\bbR^4$ intersecting in a line.}

\subsubsection{Double cover of $\sog(2,3)$ by $\spg(2,\bbR)$}\label{sec511} It was Lie who established the isomorphism between the simple Lie algebras $\soa(2,3)$ and $\spa(2,\bbR)$. This is, for example, very nicely explained in \cite{bryant}. Here we argue for this as follows:

The symplectic group $\spg(2,\bbR)$ is defined as
$$\spg(2,\bbR)=\{\glg(V)\ni A~|~\omega(Av,Aw)=\omega(v,w),\,v,w\in V\},$$
where as before $V$ is a real 4-dimensional vector space, and $\omega$ is a symplectic form on $V$. Note that $\bbZ_2=\{I,-I\}$, where $I$ is the identity in $\glg(V)$, is a subgroup of $\spg(2,\bbR)$, $\bbZ_2\subset\spg(2,\bbR)$.

Introducing, $A^\mu{}_\nu$ via $A(e_\mu)=A^\nu{}_\mu e_\nu$, and $\omega_{\mu\nu}=\omega(e_\mu,e_\nu)$, we obtain that the matrix elements of those $A\in\glg(V)$ that are in $\spg(2,\bbR)$ satisfy
\be A^\mu{}_\alpha A^\nu{}_\beta\omega_{\mu\nu}=\omega_{\alpha\beta}.\label{i1}\ee
  Since $\mathrm{dim}V=4$ we have
  $$\tfrac14\omega_{\mu\nu}\omega_{\rho\sigma}e^\mu\dz e^\nu\dz e^\rho\dz e^\sigma=\omega\dz\omega=2e^1\dz e^2\dz e^2\dz e^4=\tfrac{1}{12}\epsilon_{\mu\nu\rho\sigma}e^\mu\dz e^\nu\dz e^\rho\dz e^\sigma,$$
  and hence
\be \omega_{[\mu\nu}\omega_{\rho\sigma]}=\tfrac13\epsilon_{\mu\nu\rho\sigma}.\label{i2}\ee
  Here $\epsilon_{\mu\nu\rho\sigma}$ denotes the totally skew Levi-Civita symbol in $\bbR^4$.

  Let us now take an element $Y$ from $\omega^\perp$. We have $Y=\tfrac12 Y^{\mu\nu}e_\mu\dz e_\nu$. Then, according to \eqref{qfa} we have 
  $$-\tfrac12Q(Y)e_1\dz e_2\dz e_3 \dz e_4=Y\dz Y=\tfrac14 Y^{\mu\nu}Y^{\rho\sigma}e_\mu\dz e_\nu\dz e_\rho\dz e_\sigma=\tfrac14 Y^{\mu\nu}Y^{\rho\sigma}\epsilon_{\mu\nu\rho\sigma}e_1\dz e_2\dz e_3\dz e_4,$$ so the quadratic form $Q(Y)$ written in terms of the components $Y^{\mu\nu}=Y^{[\mu\nu]}$ of the bivector $Y$ is
  $$Q(Y)=-\tfrac12 Y^{\mu\nu}Y^{\rho\sigma}\epsilon_{\mu\nu\rho\sigma}.$$
  There is a natural action of $\spg(2,\bbR)$ on the space $\omega^\perp$ induced by the action of $\spg(2,\bbR)$ in $V$. In components it reads
  $$\spg(2,\bbR)\times\omega^\perp\quad\ni\quad (A,Y^{\mu\nu})\longrightarrow (AY)^{\mu\nu}=A^{-1}{}^\mu{}_\alpha A^{-1}{}^\nu{}_\beta Y^{\alpha\beta}\quad\in\quad\omega^\perp.$$
If we now apply the form $Q$ on the $\spg(2,\bbR)$ transformed bivector $AY$ we get
  $$\begin{aligned}
    Q(AY)=&-\tfrac12A^{-1}{}^\mu{}_\alpha A^{-1}{}^\nu{}_\beta A^{-1}{}^\rho{}_\gamma A^{-1}{}^\sigma {}_\delta Y^{\alpha\beta}Y^{\gamma\delta}\epsilon_{\mu\nu\rho\sigma}=\\&-\tfrac32A^{-1}{}^\mu{}_\alpha A^{-1}{}^\nu{}_\beta A^{-1}{}^\rho{}_\gamma A^{-1}{}^\sigma {}_\delta Y^{\alpha\beta}Y^{\gamma\delta}\omega_{[\mu\nu}\omega_{\rho\sigma]}=\\
    &-\tfrac32 Y^{\alpha\beta}Y^{\gamma\delta}\omega_{[\mu\nu}\omega_{\rho\sigma]}=-\tfrac12 Y^{\alpha\beta}Y^{\gamma\delta}\epsilon_{\mu\nu\rho\sigma}=Q(Y),
  \end{aligned}$$
  where the expressions after the second and the fourth equality sign follow from \eqref{i2}, and the expression after the third equality sign follows from \eqref{i1}. Thus the symplectic transformation $v\mapsto Av$ in $V$ induces a linear transformation $Y\mapsto AY$ in $\omega^\perp$ which \emph{preserves the real quadratic form} $Q$ of signature $(2,3)$. This gives a homomorphism of $\spg(2,\bbR)$ onto $\sog(2,3)$. Its kernel is $\bbZ_2$, since
  $$(\spg(2,\bbR)\supset\bbZ_2)\times\omega^\perp\quad\ni\quad (A=\pm I,Y^{\mu\nu})\longrightarrow(\pm\delta^\mu{}_\alpha)(\pm\delta^\nu{}_\beta)Y^{\alpha\beta}=Y^{\mu\nu}\quad\in\quad\omega^\perp.$$
  This gives the Lie's \emph{double cover} of $\sog(2,3)$ by $\spg(2,\bbR)$, $$\bbZ_2\to\spg(2,\bbR)\to\sog(2,3),$$
  which has its local version in the \emph{isomorphism of the Lie algebras} $\spa(2,\bbR)$ and $\soa(2,3)$.
  \subsection{Lie's twistor fibration}
  The relation between the groups $\spg(2,\bbR)$ and $\sog(3,2)$ recalled in the previous section is the basis for \emph{Lie's correspondence} \cite{bryant}, Section 3. This can be described in yet another incarnation of the car's fibration ${\color{red}Q}\leftarrow M\to \color{darkgreen}P$, which is \emph{Lie's twistor fibration}; see Section 4.4 in \cite{Cap} for a general theory of these things.
  
  To explain this we start with the space ${\color{red}Q}$ of all Lagrangian planes in $V$ as before. This 3-dimensional space can be locally parameterized by $(\xi,\eta,\zeta)$ as in \eqref{qc1}-\eqref{qc2}, with a Lagrangian plane ${\color{red}\dr(Y)}$ spanned by $$Y_1=(\eta+\zeta)\,e_1+e_4+\xi\,e_3\quad \&\quad Y_2=-\xi\,e_1+\,e_2+(\eta-\zeta)\,e_3.$$ There is also another 3-dimensional space associated with $V$. This is
  $${\color{darkgreen}P}:={\color{darkgreen}P(V)}=\{{\color{darkgreen}\dr(v)}~|~\lambda v,\,v\in V,\,\lambda\in\bbR\},$$
  the projectivization of $V$. This can be locally parametetrized by $(x^1,x^2,x^3)$, where a generic element of ${\color{darkgreen}P}$ is ${\color{darkgreen}\ell=\dr(x^1 e_1+x^2 e_2+x^3 e_3+e_4)}$.

  There is a third space, $M$, associated with our pair $(V,\omega)$. This is
  $$M=\{{\color{darkgreen}P}\times{\color{red}Q}\,\ni\,(\,{\color{darkgreen}\ell},{\color{red}\dr(Y)}\,)~| ~{\color{darkgreen}\ell}\in {\color{red}\dr(Y)}\},$$
  i.e. the \emph{space of all pairs} (line ${\color{darkgreen}\ell}$, Lagrangian plane associated with ${\color{red}Y}$) passing through zero in $V$ with an incidence relation such that \emph{a line} ${\color{darkgreen}\ell}$ \emph{is in the plane} ${\color{red}\dr(Y)}$. This space is \emph{four} dimensional, as a generic such pair can be parametrized by $(\xi,\eta,\zeta,s)$, where $(\xi,\eta,\zeta)$ parametrizes the plane spanned by $Y_1$ and $Y_2$, and the parameter $s$ comes from ${\color{darkgreen}\ell=\dr(Y_1+}s{\color{darkgreen} Y_2)}$ and specifies a given line from the wealth of lines passing through zero in $\color{red}\dr(Y)$.

  We again have a natural fibration ${\color{red}Q}\leftarrow M\to{\color{darkgreen}P}$:\\
  \centerline{\includegraphics[height=6cm]{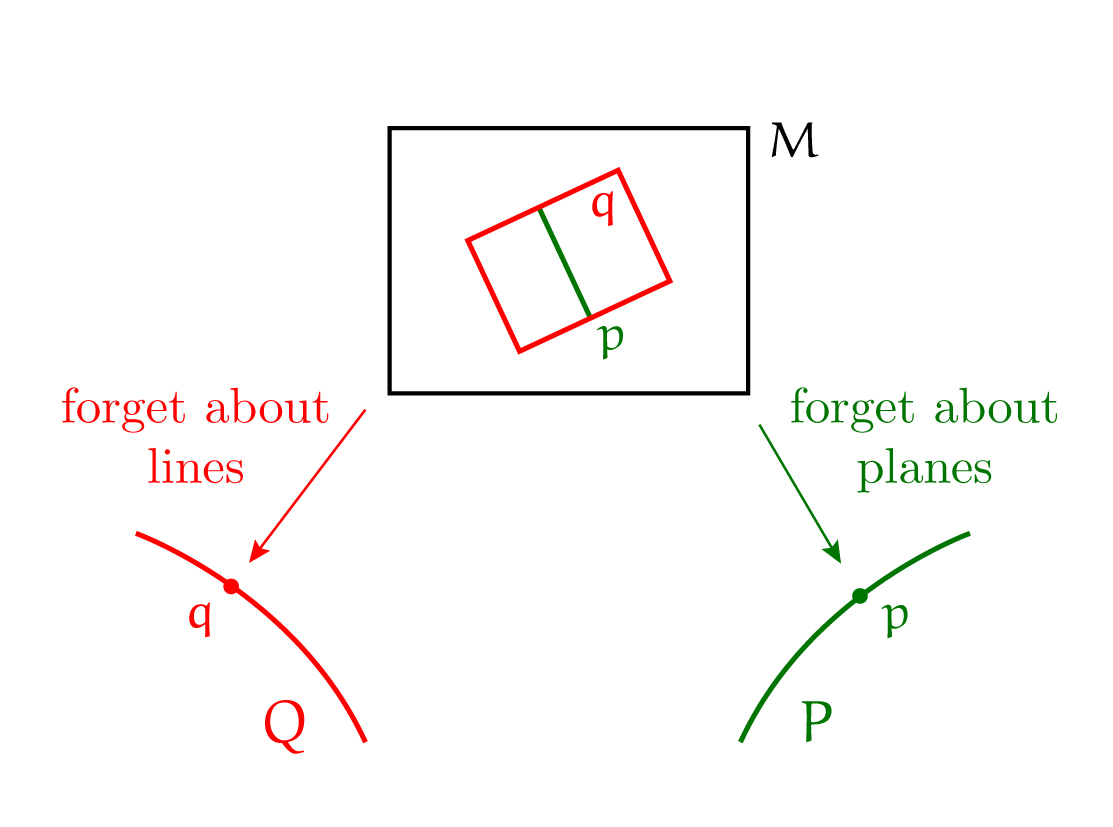}}
  where the map $M\to \color{red}Q$ is given by $({\color{darkgreen}\ell},{\color{red}Y})\to {\color{red}Y}$, and the map $M\to\color{darkgreen}P$ is given by $({\color{darkgreen}\ell},{\color{red}Y})\to {\color{darkgreen}\ell}$. It is the \emph{Lie's twistor fibration}.\\
\centerline{\includegraphics[height=6cm]{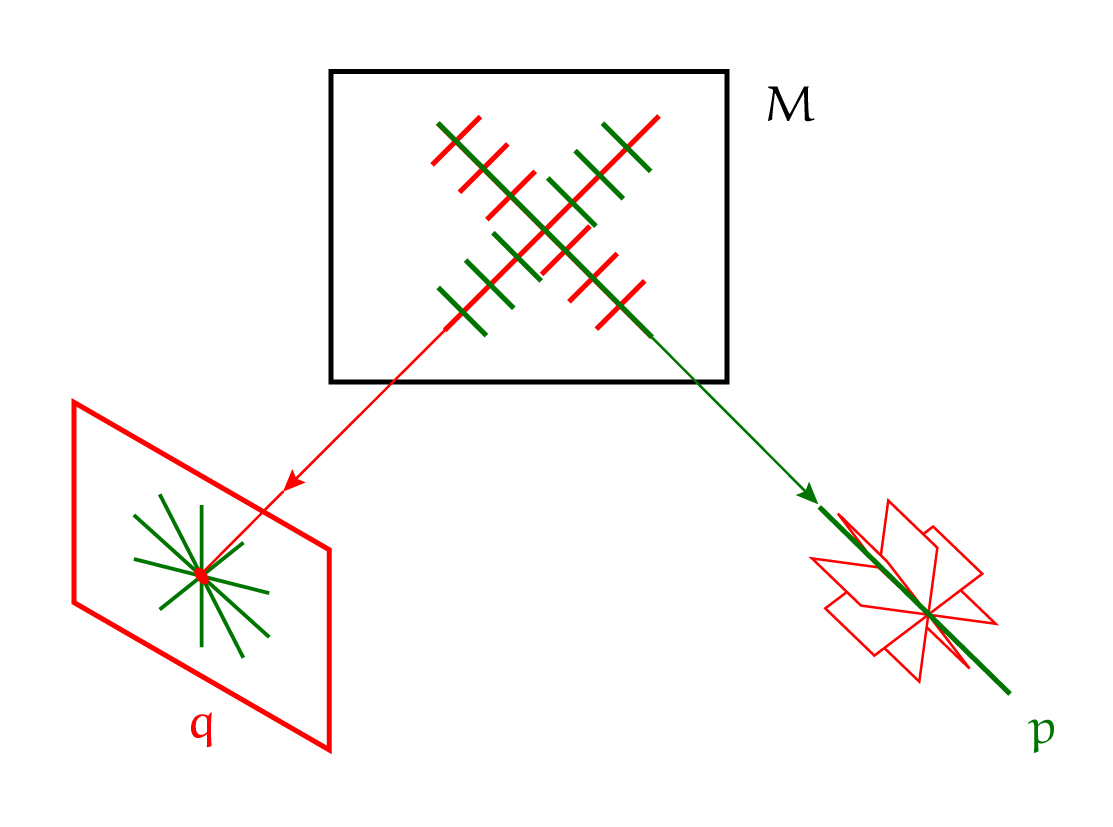}}
  In it the fiber over a point ${\color{red}q}\in \color{red}Q$, i.e. over a Lagrangian plane $\color{red}\dr(Y)$ in $V$, consists of all lines $\color{darkgreen}\ell$ passing through zero in this plane. Therefore the topology of such a fiber is the same as $\bbR P^1$. Likewise, the fiber over a point ${\color{darkgreen}p}\in\color{darkgreen}P$, i.e. over a line $\color{darkgreen}\ell$ passing through zero in $V$, consists of all Lagrangian planes $\color{red}\dr(Y)$ containing the line $\color{darkgreen}\ell$. Such a fiber also has topology of $\bbR P^1$.   

  The group $\spg(2,\bbR)$ naturally acts on ${\color{red}Q}$ and $\color{darkgreen}P$. These actions are given by
  \be
  (A,{\color{red}\dr(Y)})
  \quad\to\quad
      {\color{red}\dr(}A{\color{red}Y)}
      =
      {\color{red}\dr(\tfrac12}\,A^{-1}{}^\mu{}_\alpha A^{-1}{}^\nu{}_\beta {\color{red}Y^{\alpha\beta}\,e_\mu\dz e_\nu)}
    \label{act1}
    \ee
    and
    \be (A,{\color{darkgreen}\dr(v)})\to {\color{darkgreen}\dr(}\,A\color{darkgreen}v),\label{act2}\ee
    where $Y=\tfrac12 Y^{\mu\nu}e_\mu\dz e_\nu\in\omega^\perp$, $v\in V$ and $A\in\spg(2,\bbR)$. It also has an induced action on the elements $({\color{darkgreen}\ell},{\color{red}Y})\in M$, via
    \be (A,({\color{darkgreen}\dr(v)},{\color{red}\dr(Y)})\to ({\color{darkgreen}\dr(}A{\color{darkgreen}v)},{\color{red}\dr(}A{\color{red}Y)}).\label{act3}\ee
    It is a matter of checking that the isotropy of the action \eqref{act1} of $\spg(2,\bbR)$ on ${\color{red}Q}$ is a certain 7-dimensional group ${\color{red}P_1}$, the isotropy of the action \eqref{act2} of $\spg(2,\bbR)$ on ${\color{darkgreen}P}$ is also a certain 7-dimensional group ${\color{red}P_2}$, and that the isotropy of the action \eqref{act3} of $\spg(2,\bbR)$ on $M$ is a 6-dimensional group $P_{12}={\color{red}P_1}\cap{\color{darkgreen}P_2}$.

    Thus Lie's twistor fibration can be considered to be a double fibration of three $\spg(2,\bbR)$ \emph{homogeneous spaces}: $M=\spg(2,\bbR)/P_{12}$, ${\color{red}Q}=\spg(2,\bbR)/{\color{red}P_1}$ and ${\color{darkgreen}P}=\spg(2,\bbR)/{\color{darkgreen}P_2}$.
    \begin{align}\label{doublefib}
   \xymatrix{
        &M=\mathrm{\spg(2,\bbR)}/P_{12}  {\color{red}\ar[dl]} {\color{darkgreen}\ar[dr]} & \\
     {\color{red}Q}=\mathrm{\spg(2,\bbR)}/{\color{red}P_1} &            & {\color{darkgreen}P}=\mathrm{\spg(2,\bbR)}/{\color{darkgreen}P_2}\ . }
\end{align}   
    Due to Lie's double cover of $\sog(2,3)$ by $\spg(2,\bbR)$, and due to the proper dimensions of the spaces in the above fibration, it is clear that this gives a global version of the car's configuration space fibration
$$
    \xymatrix{
        &(M,{\color{red}{\mathcal D}}\hspace{-0.31cm}{\color{darkgreen}{\mathcal D}}={\color{darkgreen}{\mathcal D}_w}\oplus{\color{red}{\mathcal D}_g})  {\color{red}\ar[dl]} {\color{darkgreen}\ar[dr]} & \\
     {\color{red}Q} &            & {\color{darkgreen}P} . }
$$
    considered in Section \ref{dfib}. Now, the overall $\spg(2,\bbR)$ symmetry of all the ingredients of the fibration is obvious. 
    \subsection{The picture in terms of parabolic subgroups in $\spg(2,\bbR)$} The double fibration \eqref{doublefib} is a low dimensional example of the twistor correspondences discussed in \cite{Cap}, Section 4.4.6. The crucial point here is that the subgroups $\color{red}P_1$, $\color{darkgreen}P_2$ and $P_{12}$ considered in the previous section are \emph{parabolic} subgroups of a simple Lie group $\spg(2,\bbR)$; moreover they are such that $\color{red}P_1$ and $\color{darkgreen}P_2$ contain the same \emph{Borel} subgroup, which happens to be $P_{12}$. To comment about this we need some preparations.
    \subsubsection{Car's gradation in $\spa(2,\bbR)$}
    The elements $E$ of the Lie algebra $\spa(2,\bbR)$ of $\spg(2,\bbR)$ can be considered as $4\times 4$ real matrices $E=(E^\alpha{}_\beta)$ that preserve the symplectic form $\omega=\tfrac12\omega_{\mu\nu}e^\mu\dz e^\nu$, i.e. 
    $$E^\gamma{}_{\alpha}\omega_{\gamma\beta}+E^\gamma{}_{\beta}\omega_{\alpha\gamma}=0.$$
    With our choice of a basis $(e_\mu)$ in $V$, in which the symplectic form $\omega$ is as in \eqref{sympf}, the matrix $E$ giving the generic element of the Lie algebra $\spa(2,\bbR)$ is given by
    $$
    E=(E^\alpha{}_\beta)=
    \bma
    {\color{black}a_5}&{\color{darkpink}a_7}&{\color{electriccyan}a_9}&{\color{lightbrown}2a_{10}}\\
    {\color{red}-a_4}&{\color{black}a_6}&{\color{greenyellow}a_8}&{\color{electriccyan}a_9}\\
    {\color{blue}a_2}&{\color{darkgreen}a_3}&{\color{black}-a_6}&{\color{darkpink}-a_7}\\
    {\color{prune}-2a_1}&{\color{blue}a_2}&{\color{red}a_4}&{\color{black}-a_5}
      \ema,
      $$
      where the coefficients $a_I$, $I=1,2,\dots 10$, are real constants.

      Now, viewing $\spa(2,\bbR)$ as a Lie algebra consisting of all $4\times 4$ real matrices $E$ as above, with the commutator in $\spa(2,\bbR)$ being the usual commutator $[E,E']=E\cdot E'-E'\cdot E$ of two matrices $E$ and $E'$, we get a convenient basis $(E_I)$ in $\spa(2,\bbR)$ by
      $$E_I=\frac{\partial E}{\partial a_I},\quad I=1,2,\dots 10.$$
     In this basis, modulo the antisymmetry, we have the following nonvanishing commutators: $[{\color{prune}E_1},{\color{black}E_5}]=2{\color{prune}E_1}$, $[{\color{prune}E_1},{\color{darkpink}E_7}]={\color{blue}-2E_2}$,  $[{\color{prune}E_1},{\color{electriccyan}E_9}]={\color{red}-2E_4}$, $[{\color{prune}E_1},{\color{lightbrown}E_{10}}]={\color{black}4E_5}$, $[{\color{blue}E_2},{\color{red}E_4}]={\color{prune}E_1}$, $[{\color{blue}E_2},{\color{black}E_5}]={\color{blue}E_2}$, $[{\color{blue}E_2},{\color{black}E_6}]={\color{blue}E_2}$, $[{\color{blue}E_2},{\color{darkpink}E_7}]={\color{darkgreen}2E_3}$, $[{\color{blue}E_2},{\color{greenyellow}E_8}]={\color{red}E_4}$, $[{\color{blue}E_2},{\color{electriccyan}E_9}]={\color{black}-E_5-E_6}$, $[{\color{blue}E_2},{\color{lightbrown}E_{10}}]={\color{darkpink}-2E_7}$, $[{\color{darkgreen}E_3},{\color{red}E_4}]={\color{blue}-E_2}$, $[{\color{darkgreen}E_3},{\color{black}E_6}]={\color{darkgreen}2E_3}$, $[{\color{darkgreen}E_3},{\color{greenyellow}E_8}]={\color{black}-E_6}$, $[{\color{darkgreen}E_3},{\color{electriccyan}E_9}]={\color{darkpink}-E_7}$, $[{\color{red}E_4},{\color{black}E_5}]={\color{red}E_4}$, $[{\color{red}E_4},{\color{black}E_6}]={\color{red}-E_4}$, $[{\color{red}E_4},{\color{darkpink}E_7}]={\color{black}E_5-E_6}$, $[{\color{red}E_4},{\color{electriccyan}E_9}]={\color{greenyellow}-2E_8}$, $[{\color{red}E_4},{\color{lightbrown}E_{10}}]={\color{electriccyan}-2E_9}$,  $[{\color{black}E_5},{\color{darkpink}E_7}]={\color{darkpink}E_7}$, $[{\color{black}E_5},{\color{electriccyan}E_9}]={\color{electriccyan}E_9}$, $[{\color{black}E_5},{\color{lightbrown}E_{10}}]={\color{lightbrown}2E_{10}}$, $[{\color{black}E_6},{\color{darkpink}E_7}]={\color{darkpink}-E_7}$, $[{\color{black}E_6},{\color{greenyellow}E_8}]={\color{greenyellow}2E_8}$, $[{\color{black}E_6},{\color{electriccyan}E_9}]={\color{electriccyan}E_9}$, $[{\color{darkpink}E_7},{\color{greenyellow}E_8}]={\color{electriccyan}E_9}$,  $[{\color{darkpink}E_7},{\color{electriccyan}E_9}]={\color{lightbrown}E_{10}}$.
    
     What can be seen from this colorful mess?

     First, it is useful to note that our choice of the basis $E_I$ in $\spa(2,\bbR)$ is related to the following root diagram  of the Lie algebra $\spa(2,\bbR)$:\\
     \centerline{\includegraphics[height=6cm]{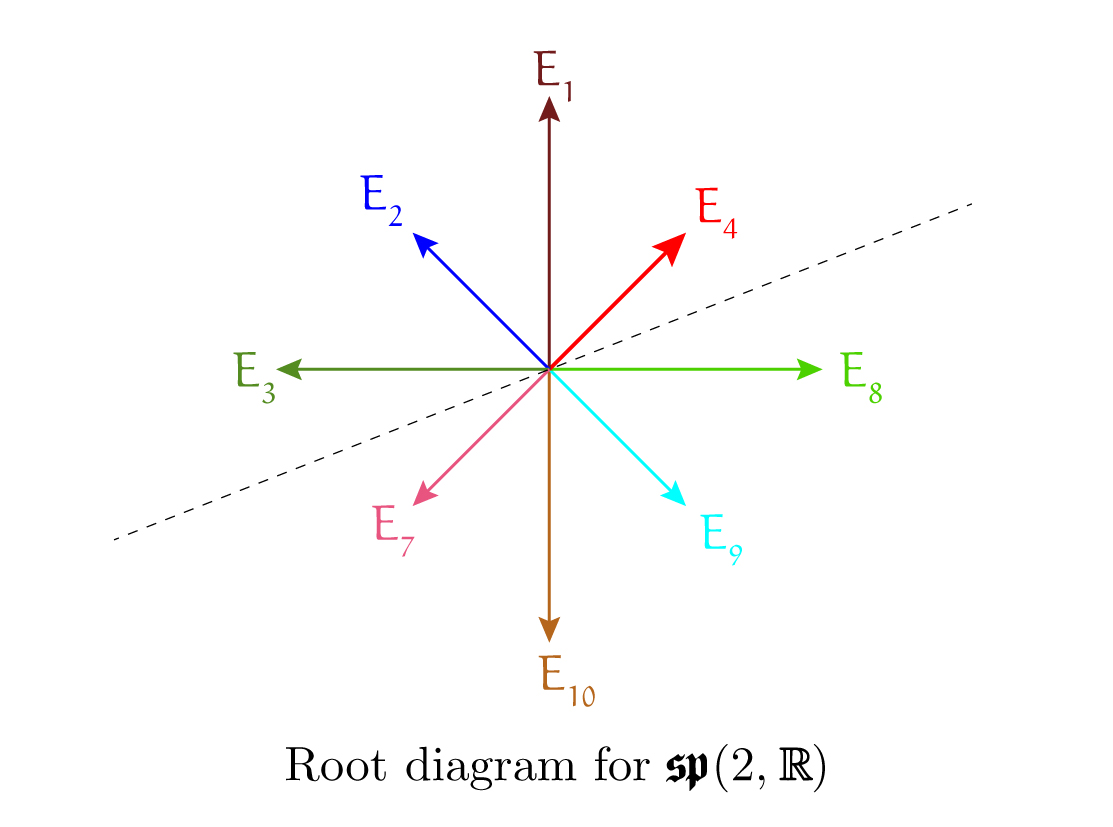}}
     This gives a mnemonic technique on how to get the directions of the vectors representing the commutators: a commutator of two vectors $E_I$ and $E_K$ in $\spa(2,\bbR)$  either \emph{vanishes} or \emph{is along the direction of} $E_I+E_K$, where the sum is the usual sum of the vectors $E_I$ and $E_K$ in the plane of the diagram. The commutators are nonzero if and only if the sum $E_I+E_K$ of vectors in the diagram belongs to the diagram.
     
     Morever, the commutation realtions above show, in particular, a certain \emph{gradation} in $\spa(2,\bbR)$.
     Indeed, define 
     $$\begin{aligned}
       {\color{prune}\mathfrak{g}_{-3}}=&\Span_\bbR({\color{prune}E_1})\\
                {\color{blue}\mathfrak{g}_{-2}}=&\Span_\bbR({\color{blue}E_2})\\
                {\color{red}\mathfrak{g}_{-1}}\hspace{-0.59cm}{\color{darkgreen}\mathfrak{g}_{-1}}=&\Span_{\bbR}({\color{darkgreen}E_3},{\color{red}E_4})\\
                {\color{black}\mathfrak{g}_{0}}=&\Span_\bbR({\color{black}E_5,E_6})\\
                {\color{darkpink}\mathfrak{g}_{1}}\hspace{-0.32cm}{\color{greenyellow}\mathfrak{g}_{1}}=&\Span_{\bbR}({\color{darkpink}E_7},{\color{greenyellow}E_8})\\
                {\color{electriccyan}\mathfrak{g}_{2}}=&\Span_\bbR({\color{electriccyan}E_9})\\
                {\color{lightbrown}\mathfrak{g}_{3}}=&\Span_\bbR({\color{lightbrown}E_{10}}),
     \end{aligned}$$
     and observe that due to the above commutation relations of the basis vectors $E_I$, these vector subspaces in $\spa(2,\bbR)$ satisfy
     $$[\mathfrak{g}_i,\mathfrak{g}_j]\subset\mathfrak{g}_{i+j},$$
     when $|i+j|\leq 3$, or
$$[\mathfrak{g}_i,\mathfrak{g}_j]=\{0\},$$
otherwise.  This observation decomposes $\spa(2,\bbR)$ onto
     $$\spa(2,\bbR)={\color{prune}\mathfrak{g}_{-3}}\oplus{\color{blue}\mathfrak{g}_{-2}}\oplus{\color{red}\mathfrak{g}_{-1}}\hspace{-0.59cm}{\color{darkgreen}\mathfrak{g}_{-1}}\oplus{\color{black}\mathfrak{g}_{0}}\oplus{\color{darkpink}\mathfrak{g}_{1}}\hspace{-0.324cm}{\color{greenyellow}\mathfrak{g}_{1}}\oplus{\color{electriccyan}\mathfrak{g}_{2}}\oplus{\color{lightbrown}\mathfrak{g}_{3}},$$
     and makes it into a \emph{3-step graded Lie algebra}.

     We further make a decomposition of ${\color{red}\mathfrak{g}_{-1}}\hspace{-0.59cm}{\color{darkgreen}\mathfrak{g}_{-1}}$ and ${\color{darkpink}\mathfrak{g}_{1}}\hspace{-0.324cm}{\color{greenyellow}\mathfrak{g}_{1}}$ onto
     $${\color{red}\mathfrak{g}_{-1}}\hspace{-0.59cm}{\color{darkgreen}\mathfrak{g}_{-1}}={\color{darkgreen}\mathfrak{g}_{-1w}}\oplus{\color{red}\mathfrak{g}_{-1g}}\quad\mathrm{and}\quad{\color{darkpink}\mathfrak{g}_{1}}\hspace{-0.324cm}{\color{greenyellow}\mathfrak{g}_{1}}={\color{darkpink}\mathfrak{g}_{1g}}\oplus{\color{greenyellow}\mathfrak{g}_{1w}}$$
     with
     $${\color{darkgreen}\mathfrak{g}_{-1w}}=\Span_\bbR({\color{darkgreen}E_3}),\,\,{\color{red}\mathfrak{g}_{-1g}}=\Span_\bbR({\color{red}E_4}),\,\,{\color{darkpink}\mathfrak{g}_{1g}}=\Span_\bbR({\color{darkpink}E_7)},\,\,\mathrm{and}\,\,{\color{greenyellow}\mathfrak{g}_{1w}}=\Span_\bbR({\color{greenyellow}E_8}).$$
     The commutation relations above show also that the following vector subspaces in $\spa(2,\bbR)$ are \emph{Lie subalgebras}:
     $$\begin{aligned}
      {\color{red}\mathfrak{p}_1}={\color{darkgreen}\mathfrak{g}_{-1w}}\oplus{\color{black}\mathfrak{g}_{0}}\oplus{\color{darkpink}\mathfrak{g}_{1}}\hspace{-0.324cm}{\color{greenyellow}\mathfrak{g}_{1}}\oplus{\color{electriccyan}\mathfrak{g}_{2}}\oplus{\color{lightbrown}\mathfrak{g}_{3}}\\
      {\color{darkgreen}\mathfrak{p}_2}={\color{red}\mathfrak{g}_{-1g}}\oplus{\color{black}\mathfrak{g}_{0}}\oplus{\color{darkpink}\mathfrak{g}_{1}}\hspace{-0.324cm}{\color{greenyellow}\mathfrak{g}_{1}}\oplus{\color{electriccyan}\mathfrak{g}_{2}}\oplus{\color{lightbrown}\mathfrak{g}_{3}}\\
       \mathfrak{p}_{12}={\color{red}\mathfrak{p}_1}\cap{\color{darkgreen}\mathfrak{p}_2}={\color{black}\mathfrak{g}_{0}}\oplus{\color{darkpink}\mathfrak{g}_{1}}\hspace{-0.324cm}{\color{greenyellow}\mathfrak{g}_{1}}\oplus{\color{electriccyan}\mathfrak{g}_{2}}\oplus{\color{lightbrown}\mathfrak{g}_{3}}\\
       \mathfrak{n}_{12}={\color{darkpink}\mathfrak{g}_{1}}\hspace{-0.324cm}{\color{greenyellow}\mathfrak{g}_{1}}\oplus{\color{electriccyan}\mathfrak{g}_{2}}\oplus{\color{lightbrown}\mathfrak{g}_{3}}\\
       \mathfrak{n}_{1}={\color{darkpink}\mathfrak{g}_{1g}}\oplus{\color{electriccyan}\mathfrak{g}_{2}}\oplus{\color{lightbrown}\mathfrak{g}_{3}}\\
       \mathfrak{n}_{2}={\color{greenyellow}\mathfrak{g}_{1w}}\oplus{\color{electriccyan}\mathfrak{g}_{2}}\oplus{\color{lightbrown}\mathfrak{g}_{3}}
     \end{aligned}$$
     \be\begin{aligned}
       \mathfrak{m}=&{\color{black}\mathfrak{g}_{-3}}\oplus{\color{blue}\mathfrak{g}_{-2}}\oplus{\color{red}\mathfrak{g}_{-1}}\hspace{-0.59cm}{\color{darkgreen}\mathfrak{g}_{-1}}\\
                {\color{red}\mathfrak{q}}=&{\color{black}\mathfrak{g}_{-3}}\oplus{\color{blue}\mathfrak{g}_{-2}}\oplus{\color{red}\mathfrak{g}_{-1g}}\\
                {\color{darkgreen}\mathfrak{p}}=&{\color{black}\mathfrak{g}_{-3}}\oplus{\color{blue}\mathfrak{g}_{-2}}\oplus{\color{darkgreen}\mathfrak{g}_{-1w}}.
                \end{aligned}\label{mqp}
     \ee
\subsubsection{Parabolic subalgebras in $\spa(2,\bbR)$}      
     We recall that a \emph{Lie subalgebra} $\mathfrak{h}$ in the Lie algebra $\mathfrak{g}$ \emph{is} ($k$-step) \emph{nilpotent} if and only if the following sequence 
     $$\mathfrak{g}_{-1}=\mathfrak{h},\quad\mathfrak{g}_{-\ell-1}=[\mathfrak{g}_{-1},\mathfrak{g}_{-\ell}],\quad\ell=1,2,\dots,$$
     of vector subspaces in $\mathfrak{g}$ terminates at step $k+1$. Here the term `terminates at step $k+1$' means that $\mathfrak{g}_{-k}\neq \{0\}$, and $\mathfrak{g}_{-k-1}=\{0\}$, for some finite $k\geq 1$. Note, that according to this definition, the \emph{Lie subalgebras} $\mathfrak{n}_{12}$, $\mathfrak{n}_1$, $\mathfrak{n}_2$, $\mathfrak{m}$, $\color{darkgreen}\mathfrak{p}$ and $\color{red}\mathfrak{q}$, of respective dimensions 4,3,3,4,3,3, \emph{are nilpotent} in $\spa(2,\bbR)$.

     Using the structure constants $c^I{}_{JK}$, defined in our basis $E_I$ of $\spa(2,\bbR)$ by $[E_I,E_J]=c^K{}_{IJ}E_K$, we find that the \emph{Killing form} $K$ of $\spa(2,\bbR)$ is
$$K=\tfrac{1}{12}K_{IJ}E^I\odot E^J=-4{\color{prune}E^1} \odot {\color{lightbrown}E^{10}}+2 {\color{blue}E^2}\odot {\color{electriccyan}E^9}+{\color{darkgreen}E^3}\odot {\color{greenyellow}E^8} -2 {\color{red}E^4}\odot {\color{darkpink}E^7}+{\color{black}E^5}\odot {\color{black}E^5}+{\color{black}E^6}\odot {\color{black}E^6},$$
     where the coefficients $K_{IJ}$ are calculated using $K_{IJ}=c^K{}_{IL}c^L{}_{JK}$. Here $E^I$, $I=1,2,\dots ,10$, is the dual basis in $\spa(2,\bbR)^*$ to the basis $E_I$ in $\spa(2,\bbR)$, $E_I\hook E^J=\delta^J{}_I$.
     
     Denoting by $\mathfrak{h}^\perp$ the subspace in $\spa(2,\bbR)$, which is \emph{Killing-form-orthogonal} to $\mathfrak{h}$,
     $$\mathfrak{h}^\perp=\{\spa(2,\bbR)\ni E~|~ K(H,E)=0,~\forall H\in\mathfrak{h}\},$$
     we can now easily see that the nilpotent subalgebras $\mathfrak{n}_1$, $\mathfrak{n}_2$ and $\mathfrak{n}_{12}$ are Killing orthogonals to the respective Lie subalgebras $\color{red}\mathfrak{p}_1$, $\color{darkgreen}\mathfrak{p}_2$ and $\mathfrak{p}_{12}$,
     $${\color{red}\mathfrak{p}_1}{}^\perp=\mathfrak{n}_1,\quad{\color{darkgreen}\mathfrak{p}_2}{}^\perp=\mathfrak{n}_2,\quad\mathrm{and}\quad \mathfrak{p}_{12}{}^\perp=\mathfrak{n}_{12}.
     $$
     Now we recall the following definition:
     \begin{definition}
       A Lie subalgebra $\mathfrak{p}$ is a \emph{parabolic} subalgebra of a (semi)simple Lie algebra $\mathfrak{g}$ if and only if its Killing orthogonal $\mathfrak{p}^\perp$ is a nilpotent subalgebra in $\mathfrak{g}$. 
     \end{definition}
     Thus according to this definition, we found \emph{three parabolic subalgebras}, $\color{red}\mathfrak{p}_1$, $\color{darkgreen}\mathfrak{p}_2$ and $\mathfrak{p}_{12}$, in the simple Lie algebra $\spa(2,\bbR)$.\footnote{It further follows that the 6-dimensional parabolic algebra $\mathfrak{p}_{12}={\color{red}\mathfrak{p}_1}\cap{\color{darkgreen}\mathfrak{p}_2}$ is a \emph{Borel} subalgebra in $\spa(2,\bbR)$.}

     \subsubsection{Twistor fibration and three flat parabolic geometries associated with a car}
     Consider now the simple Lie group $G=\spg(2,\bbR)$ and its three parabolic subgroups $\color{red}P_1$, $\color{darkgreen}P_2$ and $P_{12}={\color{red}P_1}\cap \color{darkgreen}P_2$ corresponding to the parabolic subalgebras $\color{red}p_1$, $\color{darkgreen}p_2$ and $p_{12}$ is $\spa(2,\bbR)$. Accordingly we have \emph{three} corresponding \emph{homogeneous spaces} $M=G/P_{12}$, ${\color{red}Q}=G/{\color{red}P_1}$ and ${\color{darkgreen}P}=G/{\color{darkgreen}P_{2}}$. By construction all these three spaces are $\spg(2,\bbR)$ \emph{symmetric}. Moreover, their tangent spaces at each point have the structure of the corresponding quotient vectors spaces $\mathfrak{m}=\spa(2,\bbR)/\mathfrak{p}_{12}$, ${\color{red}\mathfrak{q}}=\spa(2,\bbR)/\color{red}\mathfrak{p}_1$ and ${\color{darkgreen}\mathfrak{p}}=\spa(2,\bbR)/\color{darkgreen}\mathfrak{p}_2$. In particular, $\mathfrak{m}$, which can be identified with $\mathfrak{m}={\color{black}\mathfrak{g}_{-3}}\oplus{\color{blue}\mathfrak{g}_{-2}}\oplus{\color{red}\mathfrak{g}_{-1}}\hspace{-0.59cm}{\color{darkgreen}\mathfrak{g}_{-1}}$, has a well defined 2-dimensional vector space ${\color{red}\mathfrak{g}_{-1}}\hspace{-0.59cm}{\color{darkgreen}\mathfrak{g}_{-1}}$ with a well defined \emph{split} ${\color{red}\mathfrak{g}_{-1}}\hspace{-0.59cm}{\color{darkgreen}\mathfrak{g}_{-1}}={\color{darkgreen}\mathfrak{g}_{-1w}}\oplus{\color{red}\mathfrak{g}_{-1g}}$. This, point by point on $M=\spa(2,\bbR)$, defines an \emph{Engel distribution with a split} ${\color{red}{\mathcal D}}\hspace{-0.31cm}{\color{darkgreen}{\mathcal D}}={\color{darkgreen}{\mathcal D}_w}\oplus{\color{red}{\mathcal D}_g}$ on $M$, which \emph{by construction} is $\spg(2,\bbR)$ symmetric. Therefore this $M$ must be locally equivalent to the configuration space $M$ of a car. 

     We leave to the reader to figure out, directly from the algebraic properties of $\spa(2,\bbR)$ and $\color{red}\mathfrak{p}_1$ and $\color{darkgreen}\mathfrak{p}_2$, how the spaces ${\color{red}Q}=\spg(2,\bbR)/\color{red}P_1$ and ${\color{darkgreen}P}=\spg(2,\bbR)/\color{darkgreen}P_2$ get equipped with the respective conformal Lorentzian structure, and the contact projective structure.

     Anyhow, we can now write the global version of the car's double fibration as a \emph{parabolic twistor fibration}
     \begin{align}
   \xymatrix{
        &M=\mathrm{\spg(2,\bbR)}/P_{12}  {\color{red}\ar[dl]} {\color{darkgreen}\ar[dr]} & \\
     {\color{red}Q}=\mathrm{\spg(2,\bbR)}/{\color{red}P_1} &            & {\color{darkgreen}P}=\mathrm{\spg(2,\bbR)}/{\color{darkgreen}P_2}\  }
     \end{align}
     invoked in \eqref{doublefib}. Each space in this fibration is now a \emph{(Cartan) flat model} for a parabolic geometry of the type $(\spg(2,\bbR),P)$, where $P$ is one of $\color{red}P_1$, $\color{darkgreen}P_2$ or $P_{12}$. In this sense the car's geometry falls in the realm of \emph{parabolic geometries} \cite{Cap}.
     
     Another simpler form of this fibration, can be obtained by taking the \emph{simply connected nilpotent} Lie \emph{groups} $M$, ${\color{red}Q}$ and ${\color{darkgreen}P}$ whose corresponding Lie algebras are $\mathfrak{m}$, ${\color{red}\mathfrak{q}}$ and ${\color{darkgreen}\mathfrak{p}}$ as in \eqref{mqp}. These are \emph{Carnot groups} \cite{pansu} with additional structure, such as the Engel structure with a split on $M$. This enables us to interpret the car's double fibration as the following double \emph{fibration of Carnot groups}:
      \begin{align}
   \xymatrix{
        &M  {\color{red}\ar[dl]} {\color{darkgreen}\ar[dr]} & \\
     {\color{red}Q}&            & {\color{darkgreen}P}\ . }
     \end{align}
     Although this fibration is made only in terms of Lie groups, and although it is, in a sense, a minimal fibration locally equivalent to the car's fibration, its disadvantage with comparison to the twistor parabolic fibration is that, similar to the cars fibration, the overall $\soa(2,3)=\spa(2,\bbR)$ symmetry is not immediately visible in it.
     \subsection{Outlook: parabolic twistor fibrations in physics and in nonholonomic mechanics}
     The geometry of a car, which we discussed in this paper, is a \emph{baby version} of the well known Penrose's twistor fibration \cite{pen}
     \begin{align}
   \xymatrix{
        &M  {\color{red}\ar[dl]} {\color{darkgreen}\ar[dr]} & \\
     {\color{red}Q}&            & {\color{darkgreen}P}\ , }
     \end{align}
     in which ${\color{red}Q}$ is a 4-dimensional conformal Minkowski spacetime, ${\color{darkgreen}P}$ is a 5-dimensional space of all null rays in ${\color{red}Q}$, and $M$ is the 6-dimensional bundle of null directions over ${\color{red}Q}$ (see \cite{pen}). Penrose's fibration is also known as a basis for the \emph{Klein correspondence} (see \cite{WW}, Part 1, Section 1). To explain this we again need some preparations:

     The determination of how many different parabolic subgroups is in a given simple Lie algebra $\mathfrak{g}$, is obtained in terms of the Dynkin diagram of $\mathfrak{g}$: if $\mathfrak{g}$ is considered over the complex numbers, then the choice of a parabolic subgroup in $\mathfrak{g}$ is in one-to-one corespondence with the choice of a decoration of its Dynkin diagram with crosses marked at the nodes of the diagram. In this sense, in the Klein's-Penrose's case where the symmetry algebra is $\soa(2,4)=\sua(2,2)$, the twistor parabolic fibration looks like:\\
\centerline{\includegraphics[height=5.5cm]{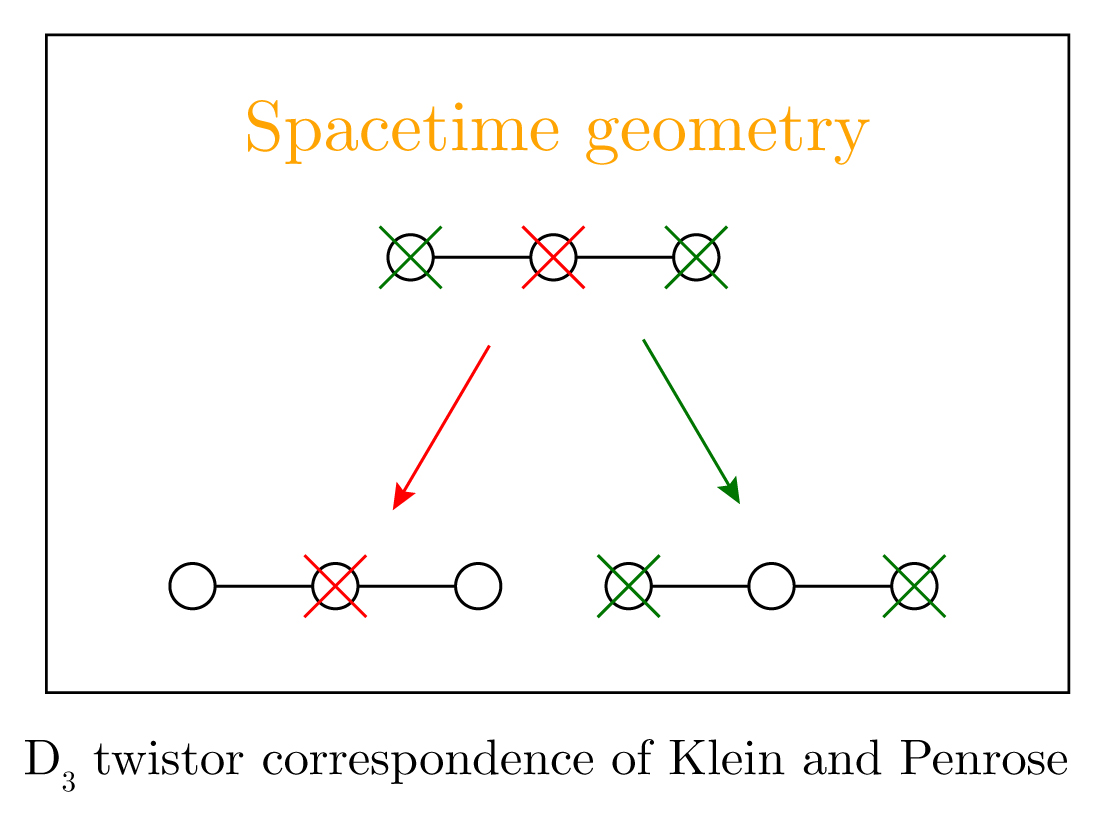}}
     Note that the symmetry Lie algebra here is a simple Lie algebra of rank 3 - there are three nodes in each of the manifolds of the diagram.

     The car's geometry is related to the symmetry algebra $\soa(2,3)=\spa(2,\bbR)$, which is a simple Lie algebra of rank 2. The corresponding twistor fibration, in terms of the Dynkin diagrams, looks like this:
     \centerline{\includegraphics[height=5.5cm]{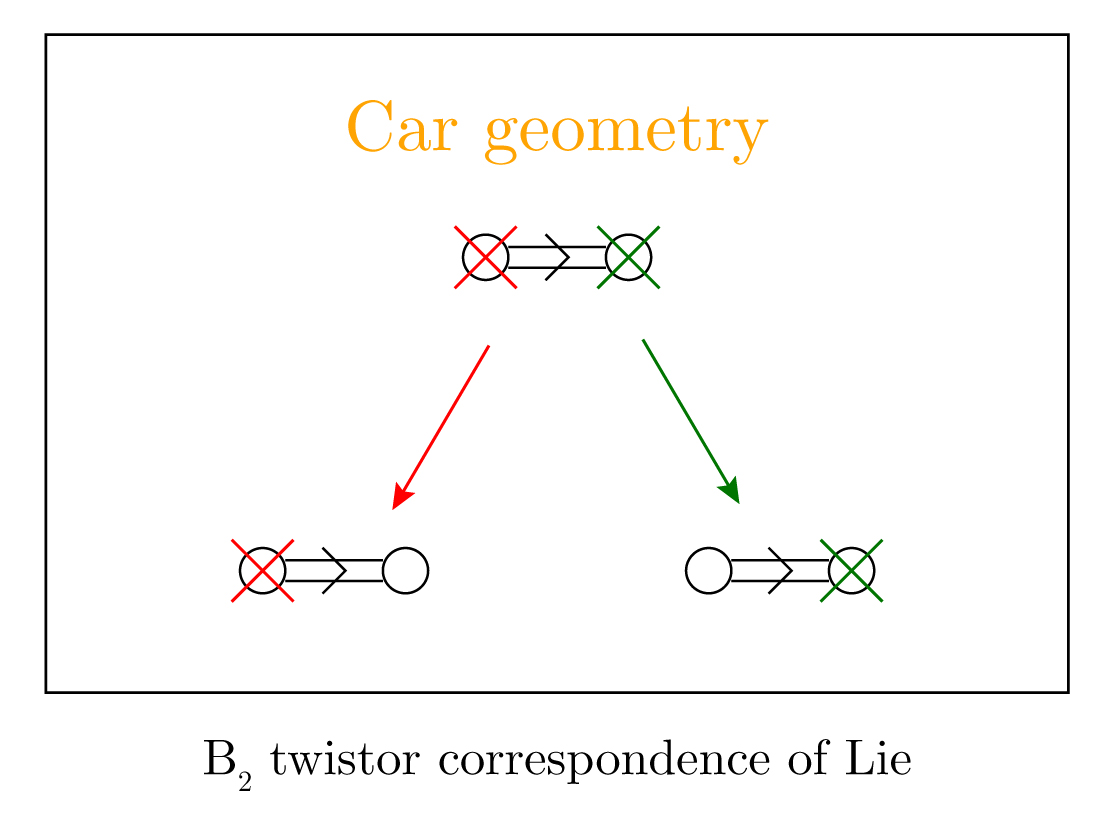}}
     R. Bryant in the beautiful article \cite{bryant} describes mathematically all the twistor parabolic fibrations associated with simple Lie algebras of rank 2. Since we interpreted the nonholonomic geometry of a car in terms of the twistor parabolic fibration related to the simple Lie algebra of Cartan-Killing type $B_2$, one can ask if there are similar physical - possibly related to nonholonomic mechanics - interpretations of the twistor parabolic fibrations related to the simple Lie algebras of type $A_2$ and $G_2$? The answer to this question is \emph{yes}.  It turns out that the $A_2$ fibration\\
     \centerline{\includegraphics[height=5.5cm]{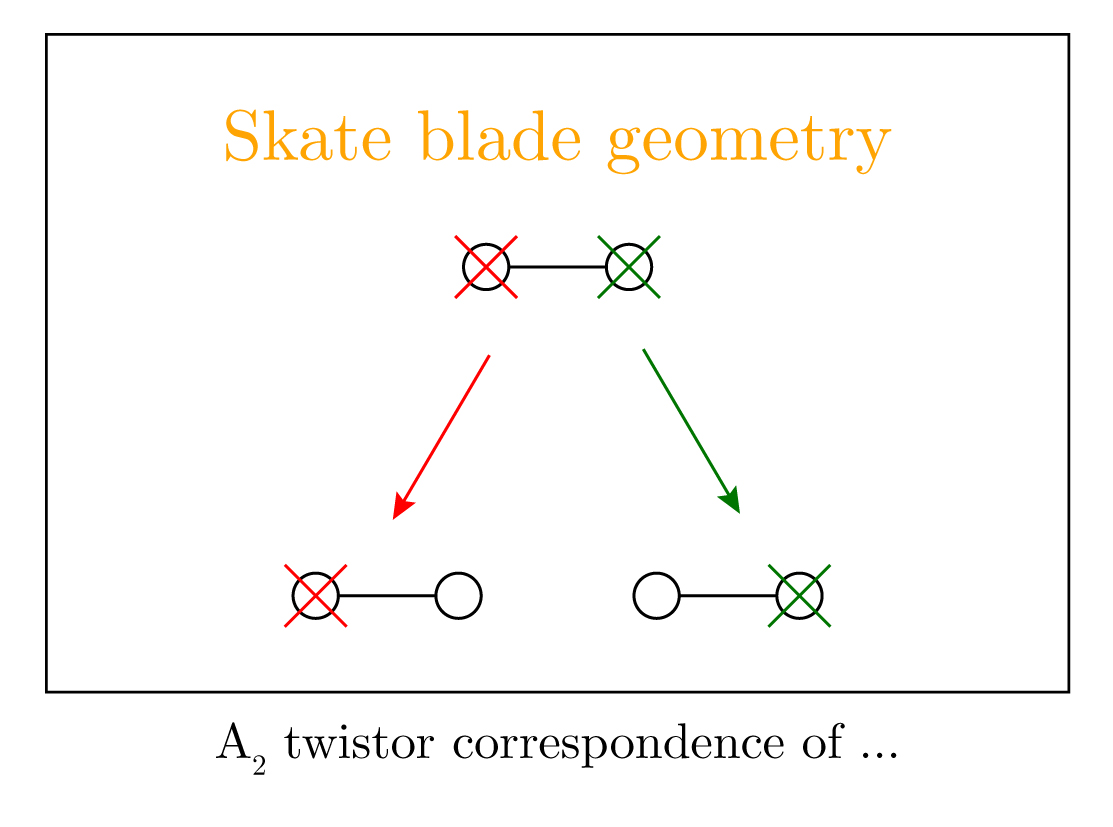}}
     corresponds to the nonholonomic movement of a \emph{skate on an ice ring} \cite{gutt,denpa}. And this case, due to the dimension of $M$ being equal to \emph{three}, and dim$\color{red}Q$=dim$\color{darkgreen}P$=2, is really the simplest to describe. It is also very similar to the car's  $B_2$ case, since $M$ is really the configuration space of the physical object (a car, a skate) subject to the nonholonomic constraints.

     The $G_2$ case \cite{cart,engel}, corresponding to the twistor diagram\\
     \centerline{\includegraphics[height=5.5cm]{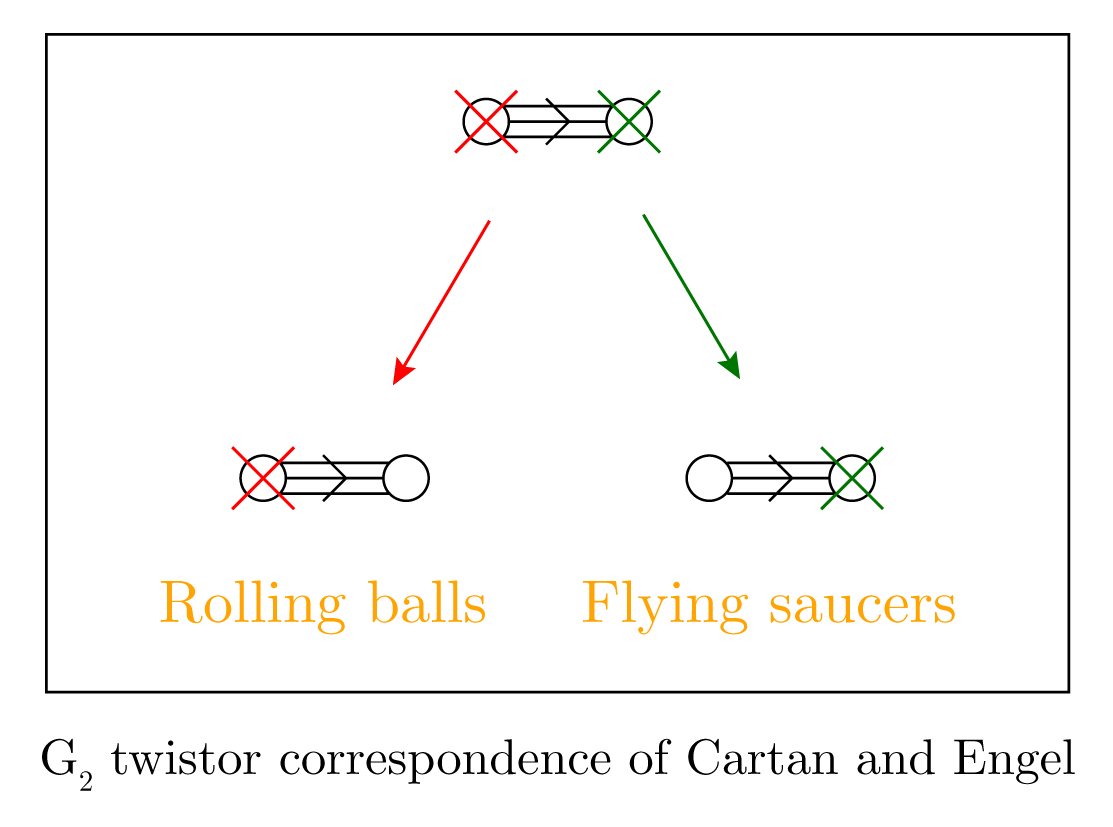}}
     is quite different. Here the dimension of $M$ is 6, and the dimensions of ${\color{red}Q}$ and ${\color{darkgreen}P}$ are both 5. In this case however, the physical objects subjected to the nonholonomic constraints (rolling surfaces, a flying saucer) have their configuration spaces as ${\color{red}Q}$ and ${\color{darkgreen}P}$ \cite{an,en1,en2} and $M$ is merely the \emph{correspondence space} enabling to translate nonholonomic movements between $\color{red}Q$ and $\color{darkgreen}P$. If in this, $G_2$ case, an interpretation of $M$ as a configuration space of some nonholonomic system exists we do not know.


\begin{thebibliography}{99}
\bibitem{an} D. An, P. Nurowski, Twistor space for rolling bodies, \emph{Commun. Math. Phys.}, {\bf 326}, 393 - 414 (2014).       
\bibitem{bryant} R. L. Bryant, \'Elie Cartan and geometric duality, in \emph{Journees Elie Cartan 1998 et 1999}, Institut Elie Cartan {\bf 16} (2000), 5--20, available at: https://services.math.duke.edu/\~{}bryant/Cartan.pdf.
\bibitem{Cap} A. \v Cap, J. Slov\'ak, (2009) \emph{Parabolic Geometries I, Background and General Theory}, Math. Surveys Monogr., vol. 154, Amer. Math. Soc.
  \bibitem{cart} \'E. Cartan, (1893) Sur la structure des groupes simples finis et continus, {\it C. R. Acad. Sci.} {\bf 116}, 784--786.
  \bibitem{car1} \'E. Cartan, (1924) Sur les varietes a connection projective, \emph{Bull. Soc. Math. France} {\bf 52} 205--241
  \bibitem{car2} \'E. Cartan (1941) La geometria de las ecuaciones diferenciales de tercer orden \emph{Rev. Mat. Hispano-Amer.} {\bf 4} 1--31
\bibitem{chern} S. S. Chern, The geometry of differential equation $y'''=F(x,y,y',y'')$, \emph{Sci. Rep. Nat. Tsing Hua Univ.}, {\bf 4} (1940), 97-111
\bibitem{en1} M. Eastwood, P. Nurowski (2019) Aerobatics of flying saucers, https://arxiv.org/abs/1810.04852
\bibitem{en2} M. Eastwood, P. Nurowski (2019) Aerodynamics of flying saucers, https://arxiv.org/abs/1810.04855
  \bibitem{engel} F. Engel, (1893) Sur un groupe simple a quatorze parametres, {\it C. R. Acad. Sci.} {\bf 116}, 786--788.
\bibitem{gstr} S. S. Chern, (1966), The geometry of G-structures,  \emph{Bull. Amer. Math. Soc.} {\bf 72}, 167-219
  \bibitem{fox} D. Fox, Contact projective structures, Indiana University Mathematics Journal 54 (2005) 1547, arXiv:math/0402332
\bibitem{godphd} M. Godli\'nski, (2008) \emph{Geometry of third order ordinary differential equations and its applications in General Relativity}, PhD thesis, Faculty of Physics, University of Warsaw  
\bibitem{gn} M. Godli\'nski, P. Nurowski, Geometry of third order ODEs, (2009) arXiv: 0902.4129,
  \bibitem{gutt} J. Gutt, Projective geometry of simple nonholonomic models, \emph{Colloquium at CFT PAS}, 24 January 2014, https://www.youtube.com/watch?v=4LV3J\_mPUAI 
  \bibitem{helgason} S. Helgason, Sophus Lie, the mathematician, In \emph{Proceedings Sophus Lie Memorial Conference, Oslo 1992}, O. A. Laudal, B. Jahren, eds. Oslo, (1994), pp. 3-21
  \bibitem{bulletin} C. D. Hill, P. Nurowski,  Differential equations and para-CR structures, \emph{Bollettino dell'Unione Matematica Italiana}, (9) III (2010) 25-91
    \bibitem{denpa} C. D. Hill, P. Nurowski, Nonholonomic mechanics and the exceptional simple Lie group $G_2$, in preparation
  \bibitem{pansu} P. Pansu, (1989), M\'etriques de Carnot-Carath\'eodory et quasiisom\'etries des espaces sym\'etriques de rang un, \emph{Annals of Mathematics}, Second Series, {\bf 129} (1): 1 -- 60.
  \bibitem{pen} R. Penrose, (1967), Twistor algebra, \emph{J. Math. Phys.} {\bf 8}, 345 -- 66
    \bibitem{WW} Ward, R., Wells, Jr, R. (1990). \emph{Twistor Geometry and Field Theory} (Cambridge Monographs on Mathematical Physics). Cambridge: Cambridge University Press. doi:10.1017/CBO9780511524493
\end{thebibliography}
\end{document}